\documentclass[11pt]{article}

\usepackage{amsmath}
\usepackage{graphicx}
\usepackage{latexsym}
\usepackage{amssymb}
\usepackage{amsbsy}
\usepackage{color}
\usepackage{multirow} 
\usepackage{rotating}
\usepackage{bm}   
\usepackage{subfigure}  
\usepackage[title]{appendix}

\usepackage{hyperref} 
\hypersetup{
    colorlinks=true,                          
    linkcolor=blue, 
    citecolor=red, 
    urlcolor=blue  } 

\newcommand{\Int}{\displaystyle\int}

\DeclareSymbolFont{msbm}{U}{msb}{m}{n}
\DeclareMathSymbol{\C}{\mathalpha}{msbm}{'103}
\DeclareMathSymbol{\R}{\mathalpha}{msbm}{'122}
\DeclareMathSymbol{\Z}{\mathalpha}{msbm}{'132}
\DeclareMathSymbol{\N}{\mathalpha}{msbm}{'116}
\newfont{\numerikNine}{cmss9}
\newfont{\numerikEight}{cmss8}

\def\RR{\mathbb R}

\def\f{\hat f}

\def\MPI{M\!P\!I}
\newcommand{\Q}{\mathcal{Q}}

\begin{document}


%
\title{ \Large
Fast Kinetic Scheme : efficient MPI parallelization strategy for 3D Boltzmann equation
\thanks{This work was supported by
the ANR JCJC project ``Moonrise'' (ANR-11-MONU-009-01)}}

%
\author{
Jacek Narski\thanks{Universit\'{e} de Toulouse; UPS, INSA, UT1, UTM;
CNRS, UMR 5219; Institut de Math\'{e}matiques de Toulouse; F-31062
Toulouse, France. ({\tt jacek.narski@math.univ-toulouse.fr})}
 }

%
\date{\today}

\maketitle

%
\begin{abstract}
  In this paper we present a parallelization strategy on distributed
  memory systems for the Fast Kinetic Scheme --- a semi-Lagrangian
  scheme developed in [J. Comput. Phys., Vol. 255, 2013, pp 680-698]
  for solving kinetic equations. The original algorithm was proposed
  for the BGK approximation of the collision kernel. In this work we
  deal with its extension to the full Boltzmann equation in six
  dimensions, where the collision operator is resolved by means of
  fast spectral method. We present close to ideal scalability of the
  proposed algorithm on tera- and peta-scale systems.
\end{abstract}

%



%
\section{Introduction} \label{sec:Intro}

Kinetic equations provide a statistical description of non equilibrium
particle gases. The evolution of the system is described by a
ballistic motion of particles interacting only by a two-body
collisions \cite{cercignani,DP-ACTA14}. The Boltzmann model derived
originally in 1870s for rarefied gases that are far from thermodynamic
equilibrium is nowadays used in variety of applications: ranging from
plasma physics to astrophysics, quantum physics, biology and social
science. In the Boltzmann description, the state of the system is
described by a distribution function defined in seven independent
dimensions: the physical space, the velocity space and the
time. Moreover, the interaction term requires multiple integrals over
velocity space to be evaluated at every space point and for every time
step of the numerical method \cite{FiMoPa:2006,PR2}. This makes the
kinetic theory very challenging from numerical view point.

There are two major strategies to
approach numerically the Boltzmann equation.
The first is to apply probabilistic methods such as Direct Simulation
Monte Carlo (DSMC) \cite{bird,Cf,CPima,Nanbu80}. The second is to
choose a deterministic scheme such as finite volume or spectral
methods \cite{DP-ACTA14,Filbet2,Filbet,Mieussens,Pal1}. The
probabilistic approach more efficient in terms of computational
time. It is however only low order with slow convergence rate.

In this work we choose a semi-Lagrangian approach
\cite{CrSon1,CrSon,Filbet2,FilbetRusso,Gu,Shoucri} applied to the
transport part of the Boltzmann equation coupled with spectral methods
to solve the collision operator
\cite{bobylev_rjasanow97,canuto:88,filbet:2011Conv,FiMoPa:2006,FilbetRusso,gamba,gamba:2010,GambaL,pareschi:1996,Villani,wu2015influence,wu2013deterministic,wu2015kinetic,wu2015fast}. In
particular, we consider a Fast Kinetic Scheme (FKS) developed originally for
the the Bhatnagar-Gross-Krook (BGK) operator \cite{BGK} in
\cite{FKS,FKS_HO,FKS_DD}. The FKS applies the Discrete Velocity Model
(DVM) technique, where the velocity space is truncated and discretized
with a set of fixed discrete velocities. As a result, the original
continuous kinetic equation is replaced by a discrete set of transport
equations that can be solved exactly in the semi-Lagrangian framework
at practically no cost. In the original method for the BGK operator,
where the collisions are modelled as a relaxation towards the local
thermodynamic equilibrium, the coupling between equations was included
in the computation of the local macroscopic variables (density,
momentum, temperature) used to approximate the local equilibrium
state. We extend herein the FKS solver to take into account more
complex collision models, such as the Boltzmann operator
\cite{bird,cercignani}. We make use of the fast spectral method
allowing to compute the collision operator in $O(N_v \log N_v)$,
where $N_v$ is a number of discrete velocity points in three
dimensions \cite{filbet:2011Conv,FiMoPa:2006,MoPa:2006}.

The curse of dimensionality makes numerical simulations of the
Boltzmann equation prohibitive on sequential machines even if fast
numerical schemes are employed. That is why the need for efficient
parallelization strategies arises. The parallel computing in the
context of kinetic equations was already explored in
\cite{Frezzotti,Frezzotti2,Frezzotti3}, where the authors made use of
Graphics Processing Unit (GPU) to solve the BGK equation with
probabilistic methods. In \cite{Malkov} the authors have implemented
the Boltzmann collision kernel on GPU, OpenMP and MPI algorithms were
explored in
\cite{Aristov,Baranger,Haack,titarev2012,titarev2014construction}. The
FKS parallelization on shared memory systems under OpenMP on GPU was
proposed in \cite{FKS_GPU} for the BGK collision kernel. In
\cite{FKS_Boltz} a simple parallelization strategy on distributed
memory systems was proposed for the FKS coupled with Boltzmann
collision kernel. The goal of this article is to study in more detail
the performance of the FKS on tera- and peta-scale systems. In
particular, we propose an efficient hybrid MPI/OpenMP implementation
of the scheme with strong scaling close to ideal on available systems.

The article is organized as follows. Section \ref{sec:FKS} introduces
the kinetic equation,
the Fast Kinetic Scheme and some collision operators.
Section \ref{sec:implementation} discusses the particle of
the FKS interpretation as particularly well suited for
parallelization, Section \ref{sec:FKS_MPI} proposes a parallel
algorithm and finally Section \ref{sec:performance} gives scalability
results and some profiling information.

\section{Kinetic equations and Fast Kinetic Scheme} \label{sec:FKS}

In the kinetic theory of rarefied gases, the state of the system is
described by a non negative distribution function $f(x,v,t)$. This
distribution function describes a density of particles moving with the
velocity $v\in \RR^{3}$ at the position $x\in\RR^{3}$ at time $t$. The
evolution of the system us governed by the six dimensional Boltzmann
equation
\begin{gather}
  \partial_t f + v \cdot \nabla_{ x }f = \Q(f,f)
  \label{eq:Boltz},
\end{gather}
where the operator $\Q (f,f)$ is the collision operator and describes
the effect of the collisions on the system. The macroscopic
characteristics (density, momentum and energy) are obtained by
integrating the distribution function multiplied by $1$, $v$ or
$|v|^2$ over the velocity space:
\begin{gather}
  U =
  \left(
  \begin{array}{c}
    \rho     \nonumber\\
    \rho u   \nonumber\\
    E
  \end{array}
  \right)
  =
  \Int_{\RR^{3}} \phi (v)f dv,
\nonumber 
\end{gather}
where the vector $\phi (v)$ is given by $(1,v,\frac{1}{2}|v|^2)^T$.  Typically,
the collision operator conserves the macroscopic quantities of the
system, {\it i.e.} the collisions preserve mass, momentum and
energy. This is expressed as
\begin{gather}
  \Int_{\RR^{3}} \phi (v)\Q(f,f) dv = 0
\nonumber 
\end{gather}
and hence the components of the vector $\phi (v)$ are referred to as {\em collision
  invariants}.

Multiplying the Boltzmann equation (\ref{eq:Boltz}) by collision
invariants and integrating over the velocity space yields a system of
equations for evolution of macroscopic conservative variables
\begin{gather}
  \frac{\partial}{\partial t}
  \Int_{\RR^{3}} f \phi (v) dv
  +
  \Int_{\RR^{3}} v\cdot \nabla_{ x }f \phi (v) dv
  =0.
  \label{eq:cons}
\end{gather}
This system is not closed as the second term involves higher order
moments. However, when the system is at thermal equilibrium, the
collision operator $\Q(f,f) =0$. The equilibrium is characterized by a
local Maxwellian distribution
\begin{gather}
  M[f] = \frac{\rho}{(2\pi T)^{3/2}}
  e^{-\frac{(u-v)^2}{2T}}
\nonumber 
,
\end{gather}
where the temperature $T$ is related to the difference between the
total and kinetic energy by the following relation:
\begin{gather}
  \frac{3}{2} \rho T = E - \frac{1}{2} \rho |u|^2
\nonumber 
.
\end{gather}
Replacing the distribution function $f$ in (\ref{eq:cons}) by the
Maxwellian distribution $M[f]$ yields a closed system --- a set of
Euler equations
\begin{align}
  &\frac{\partial}{\partial t}\rho + \nabla_x \cdot (\rho u) = 0,\nonumber\\
  &\frac{\partial (\rho u)}{\partial t}
  + \nabla_x \cdot (\rho u \otimes u+ pI) =0, \nonumber\\
  &\frac{\partial E}{\partial t}
  + \nabla_x \cdot \left( (E+p) u \right) =0,
\nonumber 
\end{align}
with the pressure following the ideal gas law $p=\rho T$.

The simplest collision operator providing the desired properties
(conservation of collision invariant and vanishing at equilibrium) is
the Bhatnagar-Gross-Krook operator \cite{BGK}
\begin{gather}
  \Q_{BGK}(f,f) = \nu (M[f] - f)
\nonumber 
,
\end{gather}
where the inter-particle collisions are modelled as a relaxation
process towards local equilibrium. The parameter $\nu = \nu (x,t)$
defines the collision frequency.

The classical Boltzmann collision operator is a multiple integral over
the whole velocity space and all possible relative angles:
\begin{gather}
  \Q_B (f,f) = \Int_{\RR^{3}} \Int_{S^2}
  B(|v-v_\star|,\theta )
  \left(
  f(v^\prime)f(v^\prime_\star) - f(v)f(v_\star)
  \right)
  d\omega dv_\star
  \label{eq:boltzmann},
\end{gather}
where $v$, $v_\star$ are velocities before collision, $v^\prime$,
$v^\prime_\star$ the velocities after collision
and $\theta $ the angle between $v-v_\star$ and
$v^\prime-v^\prime_\star$.
The post collision velocities are given by
\begin{gather}
  v^\prime = \frac{1}{2} (v + v_\star + |v - v_\star|\omega )
  \;\; , \;\;
  v^\prime_\star = \frac{1}{2} (v + v_\star - |v - v_\star|\omega ),
\nonumber 
\end{gather}
with $\omega $ being a vector on a unitary sphere $S^2$.  The
collision kernel $B$ depends only on the relative velocity before
collision and the deflection angle and has the form
\begin{gather}
  B(|v-v_\star|,\theta ) = |v-v_\star| \sigma (|v-v_\star| , \theta) 
\nonumber 
,
\end{gather}
with $\sigma $ being the scattering cross section. If the inverse $k$-th
power forces between particles are considered, $\sigma $ is given by
\begin{gather}
   \sigma (|v-v_\star|,\theta ) = b_\alpha (\theta ) |v-v_\star| ^{\alpha -1}
  \label{eq:J9hb}
\end{gather}
with $\alpha = (k-5)/(k-1)$. In the framework of the so called {\em
  variable hard spheres model} (VHS) $b_\alpha (\theta )$ is constant:
$b_\alpha (\theta )=C_\alpha $.

\subsection{Fast Kinetic Scheme}

Let us now introduce the Fast Kinetic Scheme (FKS) \cite{FKS, FKS_HO}
for solving the Boltzmann equation (\ref{eq:Boltz}). The FKS is a
semi-Lagrangian scheme \cite{CrSon1, CrSon, Filbet} which employs
Discrete Velocity Model (DVM) \cite{bobylev, Mieussens} approximation
to the original problem.

Let us start by truncating a velocity space. Let us also introduce a
cubic grid of $N_v$ equally spaced points in three dimensions. Let us
assume for simplicity that the grid step $\Delta v$ is equal in every
direction. Please note however that the FKS is not restricted to
Cartesian grids in the velocity space. In fact the method rests
unchanged even if unstructured and anisotropic velocity grids are
taken into account. Let us now define an approximation of the
continuous distribution function $f(x,v,t)$:
\begin{gather}
  \tilde{f}_k (x,t) \approx f(x,v_k,t)
\nonumber 
,
\end{gather}
that is to say, continuous $f$ is replaced by a vector $\tilde{f}$ and
the following system of $N_v$ equations is obtained:
\begin{gather}
  \partial_t \tilde{f}_k + v_k \cdot \nabla_{ x } \tilde{f}_k = \Q_k(\tilde{f},\tilde{f})
  \label{eq:DM1},
\end{gather}
where $\Q_k(\tilde{f},\tilde{f})$ is a suitable approximation of the
collision operator for the discrete velocity point $v_k$. This set of
equations is coupled only by the collision term.

Let us now discretize the physical space with a Cartesian grid
consisting of $N_s$ equally spaced points with a grid step $\Delta x$
that is equal (for simplicity reasons) in all three directions. Let us
also introduce a time discretization with $\Delta t$ being a time step
and $t^n = t^0 + n\Delta t$ for any $n\geq 0$.

In the FKS framework, the equation (\ref{eq:DM1}) is solved with a
first order splitting technique. First the transport step exactly
solves the left hand side of the problem, than the collision step
introduces the interaction using the result of the transport step as a
starting point:
\begin{eqnarray}
  \label{eq:transport} \textit{Transport stage}  & \longrightarrow & \partial_t \tilde{f}_{k} + v_{k} \cdot\nabla_{x} \tilde{f}_{k} = 0, \\
\nonumber 
 \textit{Collisions stage} & \longrightarrow & \partial_t \tilde{f}_{k}  = \Q_k(\tilde{f},\tilde{f}).
\end{eqnarray}
Please note that higher order splitting techniques may also be
considered.

\subsubsection*{Transport step}

Let $f^n_{j,k}$ be a point-wise approximation of the distribution
function at time $t^n$, position $x_j$ and velocity $v_k$: $f^n_{j,k}
= f(x_j,v_k,t^n)$. The main idea behind FKS is to solve the transport
step (\ref{eq:transport}) exactly. Let us define a piece wise constant
in space approximation $\bar{f}^n_k(x)$ of the function
$\tilde{f}_k(x,t^n)$ such that $\bar{f}^n_k(x) = f(x_j,v_k,t^n)$ if
$x\in [x_{j-1/2} , x_{j+1/2}]=\Omega _j$ belongs to the space cell
centered on $x_j$. The exact solution to the transport step at time
$t^n$ is therefore given by
\begin{gather}
  \bar{f}_k^{\star,n+1} = \bar{f}^n_k(x - v_k \Delta t)
\nonumber 
.
\end{gather}
The function $\bar{f}^n_k$ is advected with a velocity $v_k$ during a
time step $\Delta t$.
The
discontinuities of $\bar{f}^{\star,n+1}_k$ does coincide now with
space cell boundaries after the transport step.

\subsubsection*{Collision step}

During the collision step the amplitude of the distribution function
$\bar{f}$ is modified. The collision operator is solved locally on the
space grid points and than extended to the whole domain $\Omega $. The
following equations (ordinary differential or integro-differential)
are solved:
\begin{gather}
  \partial_t f_{j,k} = \Q_k(f_{j,\cdot},f_{j,\cdot})
\nonumber 
,
\end{gather}
where $f_{j,k} = f(x_j,v_k,t)$ for all space and velocity grid points
$j=1,\ldots ,N_s$ and $k=1,\ldots , N_v$ and $f_{j,\cdot}$ is a vector
representing the distribution function at the space cell $j$ composed
of $f_{j,k}$. The initial data for this equation are provided by the
transport step performed before. The time discretization chosen in
this work is the first order explicit Euler scheme
\begin{gather}
  f^{n+1}_{j,k} = f_{j,k}^{\star,n+1} + \Delta t \Q_k(f_{j,\cdot}^{\star,n+1},f_{j,\cdot}^{\star,n+1})
  \label{eq:disc_relax},
\end{gather}
where $f_{j,k}^{\star,n+1} = \bar{f}_k^{\star,n+1}(x_j)$ is the value
of transported distribution function at grid point $x_j$ and
$f_{j,\cdot}^{\star,n+1}$ is a vector composed of
$f_{j,k}^{\star,n+1}$.  Please note that other type of time
integrators can be successfully implemented instead of this forward
scheme. In particular the special care must be taken in the stiff
limit, please refer to \cite{Dimarco_stiff1,Dimarco_stiff2}.

Equation (\ref{eq:disc_relax}) furnishes a modified value of the
distribution function at grid points $x_j$ for velocity points $v_k$
at time $t^{n+1}$. In order to obtain the value of $f$ at every point
of the domain a new piecewise constant function $\bar{\Q}_k$ is
defined for every discrete velocity $v_k$:
\begin{gather}
  \bar{\Q}_k^{n+1} (x) =
  \Q_k(f_{j,\cdot}^{\star,n+1},f_{j,\cdot}^{\star,n+1})
  \;\; , \;\;  \forall x \;\; \text{such that} \;\;
  \bar{f}_k^{\star,n+1}(x) = f_{j,k}^{\star,n+1}
\nonumber 
,
\end{gather}
that is to say, the collision operator at every point of $\Omega $ is
approximated by a piecewise constant function with discontinuities
located at the same points as the piecewise constant function that
approximates the distribution function after the transport
step. Thanks to this assumption, the spatially reconstructed
distribution function after the collision step reads
\begin{gather}
  \bar{f}^{n+1}_k (x) = \bar{f}_k^{\star,n+1}(x) + \Delta t \bar{\Q}_k^{n+1} (x)
\nonumber 
.
\end{gather}
This completes the description of the Fast Kinetic Scheme.

\subsection{Collision operator}

Let us now focus on the details related to the collision kernel.

\subsubsection{BGK approximation}

If particle interaction is modeled by relaxation towards local
equilibrium, the collision term $\Q_k
(f_{j,\cdot}^{\star,n+1},f_{j,\cdot}^{\star,n+1})$ becomes $\nu ({\cal
  E}_{j,k} - f_{j,k}) $, where ${\cal E}_{j,k}$ is a suitable
approximation of the Maxwell distribution for the velocity $v_k$ at
the grid point $x_j$. As the Maxwellian distribution depends on the
macroscopic characteristics of the system that are unchanged during
the relaxation step (since they are {\em collision invariants}), the
relaxation step (\ref{eq:disc_relax}) becomes completely decoupled. In
particular, $\Q_k$ depends only on one velocity point $v_k$ and not on
the others.

\subsubsection{Boltzmann operator}

If the Boltzmann operator is considered, the collision operator $\Q_k
(f_{j,\cdot}^{\star,n+1},f_{j,\cdot}^{\star,n+1})$ involves
integration over whole velocity space for every point $x_j$ of the
space grid. The relaxation step is solved by means of Fast Spectral
Scheme presented in the Appendix \ref{Boltz} and requires multiple
Fourier transforms to be computed at every time step and at every
space cell.

\section{Implementation}\label{sec:implementation}

Let us switch to a particle interpretation of the FKS. Every point of
the velocity grid represents a particle moving with velocity
$v_k$. Every space cell $\Omega _j$ centered on the space grid point
$x_j$ contains exactly the same set of particles at exactly the same
relative positions. Therefore one needs to store the particle position
and velocity only in one generic cell and not in the whole
domain. This reduces the memory requirements seven times: only mass of
the particles is stored for every point of the $6D$ grid, three
components of particle position and velocity vectors are only required
for the generic reference cell. The distribution function is related
to particle masses by
\begin{gather}
  f(x,v,t) = \sum_{j,k=1}^{N_s,N_v} \mathfrak{m}_{j,k}(t) \, \delta(x-x_{j,k}(t))\delta(v-v_{j,k}(t)), \quad v_{j,k}(t)=v_k,
\nonumber 
\end{gather}
where $x_{j,k}(t)$ is particle position, $v_{j,k}(t)$ is its velocity
and $\mathfrak{m}_{j,k}(t)$ particle mass. In the FKS the particle
velocity is unchanged and the position is altered during the transport
step:
\begin{gather}
  x_{j,k}(t+\Delta t)=x_{j,k}(t)+v_{j,k}(t)\Delta t.
\nonumber 
\end{gather}
The transport step moves the particles in the reference cell. The
motion of particles in the remaining cells is identical. If a given
particle escapes the generic cell, another one with the same velocity
(but different mass) enters the cell from the opposite side.

The collision step modifies the particle masses in every space cell:
\begin{gather}
  \mathfrak{m}_{j,k}(t+\Delta t)= \mathfrak{m}_{j,k}(t) \, + \, \Delta t \, Q_k(v_{j,\cdot}),
\nonumber 
\end{gather}
where $Q_k(v_{j,\cdot})$ is the approximation of the collision
operator in the center of the cell by the means of the fast spectral
method presented above.

The macroscopic variables at time $t^n$ are defined on the space grid only and are
computed as a sum over particles in the given cell $j$:
\begin{gather}
  U_j^n
  =\sum_{k=0}^{N_v}
  \phi (v_{j,k})
  \mathfrak{m}_{j,k}^n
  (\Delta v)^3
\nonumber 
.
\end{gather}
As the collision step does not change the macroscopic conservative
variables, they can be efficiently computed at time $t^{n+1}$
after the transport step by adding the contribution from the particles
leaving and entering the given cell $j$ to the values at the previous
time step. If a particle $(j,k)$ is transported to the cell $j+\delta$
during the transport step, there is a sister particle entering the
cell $j$ from $j-\delta $. A suitable contribution has to be added to
from the conservative variables in the cell $j$ :
\begin{gather}
  U_j^{n+1}
  =
  U_j^n
  +
  \sum_{k, \ x_{j,k}^{n+1} \in \Omega_{j+\delta}, \ x_{j,k}^{n} \in \Omega_{j}}
  (\mathfrak{m}_{j-\delta ,k}^n - \mathfrak{m}_{j,k}^n )
  \phi (v_{j,k})
  (\Delta v)^3
  \label{eq:DM_update}
\end{gather}

The most expensive part in the algorithm is the collision
operator. Even in the case of the relatively simple BGK approximation
the computation of the relaxation term takes $90\%$ of the
computational time on serial machines \cite{FKS_GPU}. The cost of the
Boltzmann integral is substantially greater, even if the Fast Spectral
Method is employed. Indeed, a number of FFTs must be performed for
every space cell and for every discrete angle in order to compute
convolutions. If 16 discrete angles are considered ($A_1 = A_2 =4$),
this number equals 96 and even if the size of those transforms is
relatively small, this represents the main computational burden.
Evaluation of the Boltzmann operator represents more than $99\%$ of
the computational time on serial machines. Fortunately the collision
operator is in some sense decoupled from the FKS framework: it can be
implemented independently of the FKS. This suggest a following
strategy for the parallelization on distributed memory systems. On the
upper level, the FKS is parallelized with MPI over available
computational nodes. On the collision level, a suitable operator is
implemented on the available node architecture: using the classical
OpenMP type parallelism or the SIMD (Single Instruction, Multiple
Data) programming model on GPUs or on the Intel Many Integrated Core
(MIC) hardware. The implementation details of the collision operator
does not influence the MPI scalability of the algorithm.

\section{Fast Kinetic Scheme on parallel machines} \label{sec:FKS_MPI}

In this section we propose a simple parallelization strategy for
distributed memory systems. There are two possible approaches: the
first one is to decompose velocity space and distribute it over
computational nodes keeping all spatial degrees of freedom at every
node. This approach is very similar to the strategy employed in the
GPU algorithm for the BGK collision kernel ({\it cf.}  \cite{FKS_GPU}
), where for each velocity point from the velocity mesh a relaxation
term was computed in parallel for all space cells at once on a GPU
device. It was also chosen in \cite{titarev2012} for the MPI
implementation.  Every computational node performs computations of a
relaxation for a subset of velocity grid. Then the partial moments are
evaluated. The total moments are obtained from gathering all
contributions from all computational nodes. This approach is well
suited for collision kernels that are local, {\it e.g.} for the BGK
approximation, where the collision computed for a given position in
physical and velocity space depends only on the distribution function
at the same position and on total moments. As every node contains all
spatial degrees of freedom, no particle escapes given computational
node and no particle mass is exchanged with neighbouring nodes. The
MPI communication is limited only to the partial moments. Another
advantage is that even if a complicated domain is considered {\it
  i.e.}  containing perforations, no complicated domain decomposition
or load balancing techniques are required to ensure equal workloads
across computational nodes. However, any collision kernel that is
non-local would generate huge amount of communication between all
processors. Boltzmann collision operator involve a double integral
over velocity space meaning that at every iteration every node must
have access to the whole velocity space. This approach is therefore
not well suited for the Boltzmann operator.

The second possibility, adopted herein, is to distribute spatial
degrees of freedom over computational nodes, keeping on every mode a
complete velocity space.  Since the update of conservative variables
(density, momentum and energy) requires an exchange of particle mass
with neighbouring spatial cells (and does not involve any summation
over whole physical space), the internodal memory transfer is limited
to merely cells located on a boundary of a subdomain.  Moreover, the
information is exchanged with one node only and not with every node
reserved for the computation. Comparison of the two approaches can be
found in \cite{titarev2014construction}.

The spatial domain is decomposed into equally sized non-overlapping
cuboids, pencils or slabs with ghost layers. Depending on the choice,
every node has to communicate with 2 (for slabs), 8 (for pencils) or
26 (for cuboids) neighboring nodes. Cuboids usually minimize the size
of ghost layers but have the biggest MPI overhead as they require more
calls to MPI in order to communicate with all neighbors. Thus, the
cuboid domain decomposition strategy is not necessarily optimal and
better results can be sometimes obtained when using pencils or slabs.

Moreover, the performance of the method can be improved by OpenMP or
SIMD parallelization applied to the loops over velocity space at each
node. This kind of parallelization for shared memory systems was
already proposed in \cite{FKS_GPU}. Let us now, for the sake of
completeness, repeat the sequential and parallel algorithms for shared
memory systems.

\subsection{Sequential algorithms for FKS}

We consider particle positions $\bm{X}_p^{n}$ and masses $\mathfrak{m}_{j,p}^n$
known at time $t^{n}$ as well as conservative $F_j^n$ and primitives
variables $(\rho, \bm{U}, T)_j^{n}$. The algorithm
reads:
\begin{enumerate}

\item \textit{Transport of particles.} Displace $N_v$ particles with
  (\ref{eq:transport}), produce a list of $N_{\text{out}}$ particles
  escaping the generic cell and store the $\delta$s determining the
  destination and provenance of associated sister particles.

\item \textit{Update conservative variables} $U_j^{n+1}$ using
  (\ref{eq:DM_update}) and the results from the transport step.

\item \textit{Relaxation step.} Compute masses of $N_v$ particles with
  a collision kernel of choice, store them in an array of the size
  $N_v \times N_s$.

\end{enumerate}

\subsection{Classical parallel architecture: Open-MP} \label{ssec:open-MP}

The modified algorithm reads:
\begin{enumerate}

\item \textit{Relaxation step.} Divide the number of spatial cells by
  the number of processors. Compute in parallel the masses of $N_v$
  particles with a collision kernel of choice, parallelization is
  performed on the loop over the number of mesh points in the physical
  space. This computation is local on the space mesh.

\item \textit{Transport of particles.} Move in parallel $N_v$
  particles with (\ref{eq:transport}). This step is done in only one
  space cell. The motion of particles in the other cells is the same.

\item \textit{Update conservative variables.} Test in a parallel loop
  over the number of mesh points in the physical space if a particle
  has escaped from the generic cell. If so, add a contribution to
  $U_j^{n+1}$ using (\ref{eq:DM_update}) for every space cell. Update
  the particle position and exchange particle mass with the associated
  sister particle.
  
\end{enumerate}

\subsection{Graphic Processing Unit (GPU) architecture: CUDA} \label{ssec:open-GPU}

This parallelization design can be summarized in the following
algorithm.
\begin{enumerate}
\item \textit{Copy from CPU to GPU}. Copy to the GPU memory all
  primitive and conservative variables.

\item \textit{Loop over $N_v$ particles}

  \begin{enumerate}
  
  \item \textit{Relaxation step} Compute relaxed masses of particles
      for every space cell using CUDA. Store the result on GPU.
      
    \item \textit{Transport step} Move every $N_v$ particle and test if it
      has escaped the generic cell. If so, store the provenance cell of the
      sister particle.

    \item \textit{Update conservative variables}. If the particle has
      escaped the generic cell, add contribution to conservative
      variables. Reassign its mass and position with the ones of the
      incoming sister particle.
      
    \item \textit{Copy from GPU to CPU.} Copy the resulting mass array
      from the GPU memory to the CPU memory.
      
  \end{enumerate}
  
\item \textit{Copy from GPU to CPU}. Write to the CPU memory the
  updated conservative and primitive variables.

\end{enumerate}

\subsection{MPI version of FKS}

The parallel algorithm is straightforward
\begin{enumerate}

\item \textit{Initialization} Divide the computational domain into
  $N_{\MPI} =N_{\MPI_x}\times N_{\MPI_y} \times N_{\MPI_z}$ equally
  sized cuboids. Allocate memory on each computational node: arrays of
  the size $N_s/ N_{\MPI} $ for storing the conservative and primitive
  variables relative to a given subdomain and an array of the size
  $N_v \times (N_s/N_{\MPI} + N_{\text{ghosts}} )$ for storing masses
  relative to a given subdomain with additional ghost layers
  containing masses of particles in the adjacent space cells.

\item \textit{Time iterations} For every computational node:

  \begin{enumerate}

  \item \textit{Relaxation step} performed in parallel on GPU or with
    OpenMP for every particle in a given subdomain.

  \item \textit{Transport of particles.} Move in parallel $N_v$
    particles with (\ref{eq:transport}). This step is done in only one
    space cell in each subdomain. The motion of particles in the other
    cells is the same.
    
  \item \textit{Communication} If a particle is escaping from a given
    subdomain, broadcast its mass to suitable computational node.
    
  \item \textit{Update conservative variables}. If the particle has
    escaped the generic cell, add contribution to conservative
    variables. Reassign its mass and position with the ones of the
    incoming sister particle. For particles located on the boundary of
    the cuboid and escaping the subdomain use the values stored in the
    ghost cells in the previous step.

  \end{enumerate}

\item \textit{Finalization} Free memory and close MPI communication.

\end{enumerate}

The communication is limited only to neighboring subdomains.  The
amount of data to be exchanged depends only on the local mesh size and
chosen MPI topology (slabs~/ pencils~/ cuboids). It does not depend
on number of computational nodes employed. Weak scaling is therefore
evident.

\section{Performance tests} \label{sec:performance}

Numerical tests for the BGK collision kernel were performed on the EOS
supercomputer at CALMIP, Toulouse. The supercomputer is equipped with
612 computational nodes, each of them containing two
Intel$^\text{\textregistered}$ Ivy Bridge $2.8$GHz 10 core CPUs and 64
GB of RAM. Each CPU was equipped with $25$MB of cache memory. The code
was compiled with gcc-5.3.0 and executed on 2 to 90 computational
nodes. That is to say, on 40 to 1800 computational cores in parallel.
The tests for Boltzmann collision kernel were performed on the EOS
supercomputer ($N=64^3$ meshes) and on the thin nodes of GENCI-TGCC
supercomputer Curie for $N=128^3$. The machine is equipped with 5040
B510 Bullx nodes (called thin nodes), each containing two
Intel$^\text{\textregistered}$ Sandy Bridge $2.7$GHz 8 core CPUs (20MB
of cache memory) and 64GB of RAM. In the case of the Boltzmann
collision operator the Fast Fourier Transforms were computed by means
of the fftw library, version 3.3.4. The code was executed with 20
OpenMP threads per node on EOS and 16 OpenMP threads per node on
Curie.

The performance of the parallel algorithm was tested on the 3D Sod
test case. The problem description is a 3D explosion problem
\cite{toro-book}, the initial state being given by the well known Sod
shock tube problem.  Let us consider a cubic domain of size $[0,2]^3$.
Left and right states of the 1D Sod problem are given by a density
$\rho_{L}=1$, mean velocity $\bm{U}_{L}=\bm{0}$ and temperature
$\theta_{L}=5$, while $\rho_{R}=0.125$, $\bm{U}_{R}=\bm{0}$,
$\theta_{R}=4$.  The gas is initially in thermodynamic equilibrium.
The left state is set for any cell inside a ball centered in $(1,1,1)$
and of radius $0.2$. The right state is set elsewhere. The
computations are stopped at final time $t_\text{final}= 0.07$ for the
BGK operator and after $16$ iterations for the Boltzmann
operator. This number of times steps ensured that every discrete
particle has changed the physical cell at least once. We consider the
case in which $\tau =10^{-1}$, {\it i.e.} far from the fluid limit. We
are not interested in the convergence of the numerical solution but in
the parallel efficiency only. Some convergence results and more
interesting (from physical point of view) test cases are presented in
\cite{FKS_Boltz}.

The performance tests were run for the BGK and 3D Boltzmann collision
kernels. Computation of the 3D Boltzmann collision operator is much
more time consuming compared to the BGK operator. The run time is 100
times bigger for $N_v=32^3$ velocity points. The ratio increases even
more for larger velocity meshes. We expect that all the communication
time and MPI overhead will be hidden for the 3D Boltzmann operator
with MPI efficiency close to one on the supercomputer at our
disposal. When dealing with the BGK kernel it is the memory
requirement that is a real bottleneck. Indeed, BGK $3D\times 3D$
simulations on medium size meshes can be efficiently performed on
multi-core or GPU based shared memory systems with enough RAM at
disposition with runtime that of the order of hours or days rather
than weeks or months \cite{FKS_GPU}. However, if fine scale
simulations are required the MPI parallelization is indispensable.

\subsection{BGK}

\begin{figure}[!t] 
  \centerline{\resizebox{0.4\textheight}{!}{\rotatebox{0}
      {\includegraphics{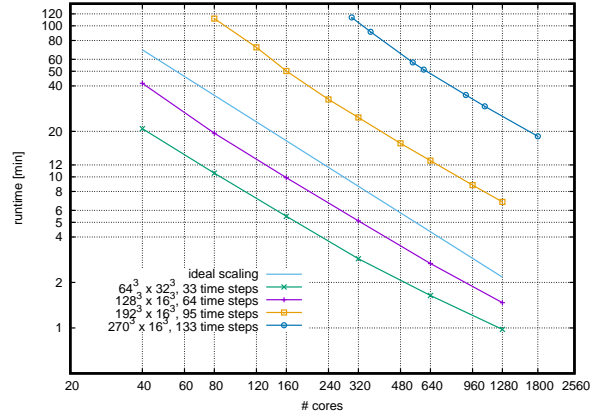}}}}
  \caption{Computational time as a function of number of cores employed for BGK.}
  \label{fig:cores_runtime}
\end{figure}
\begin{figure}[!t] 
  \centerline{
    \resizebox{0.4\textheight}{!}{\rotatebox{0}
      {\includegraphics{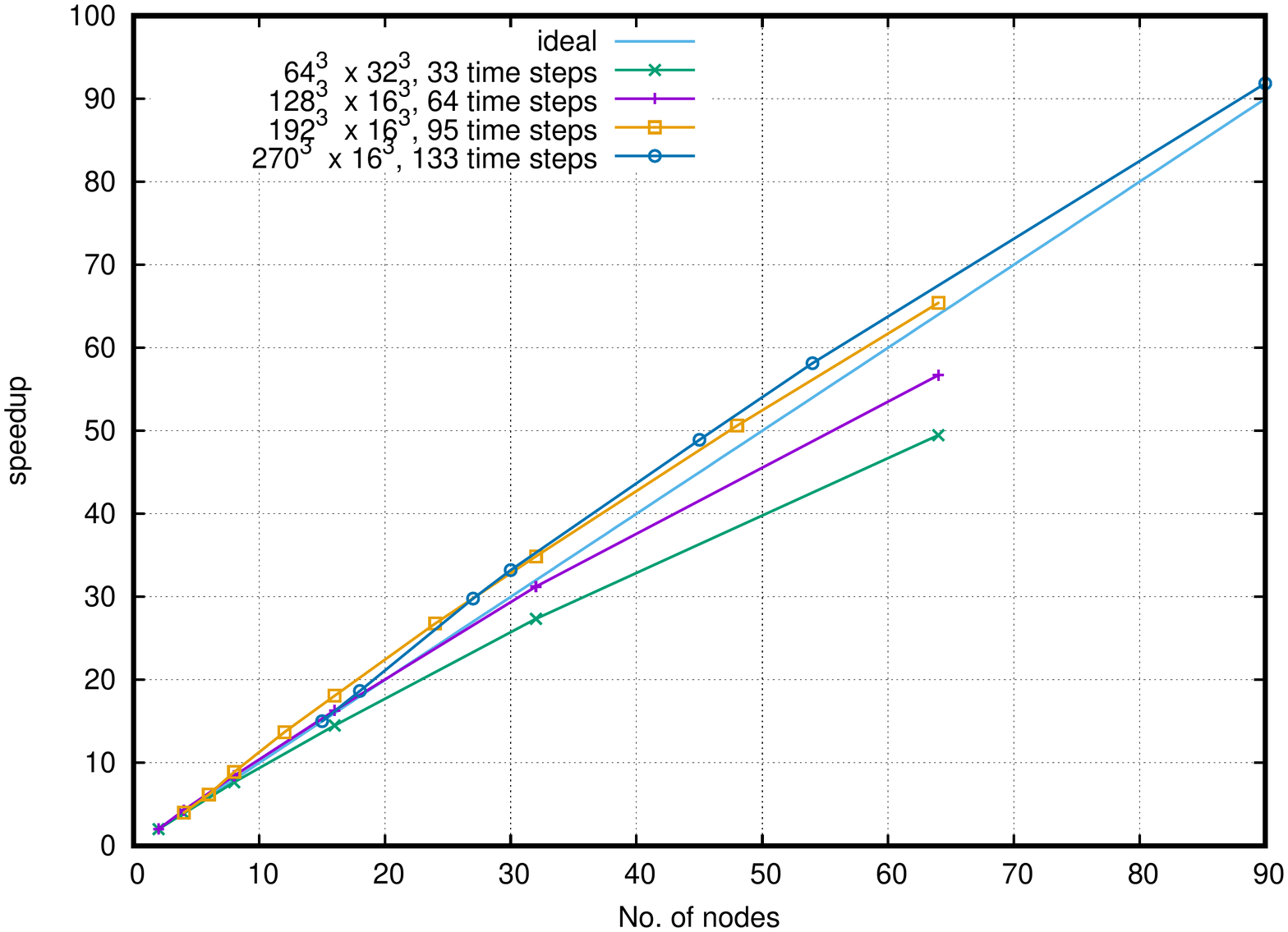}}}
    \resizebox{0.4\textheight}{!}{\rotatebox{0}
      {\includegraphics{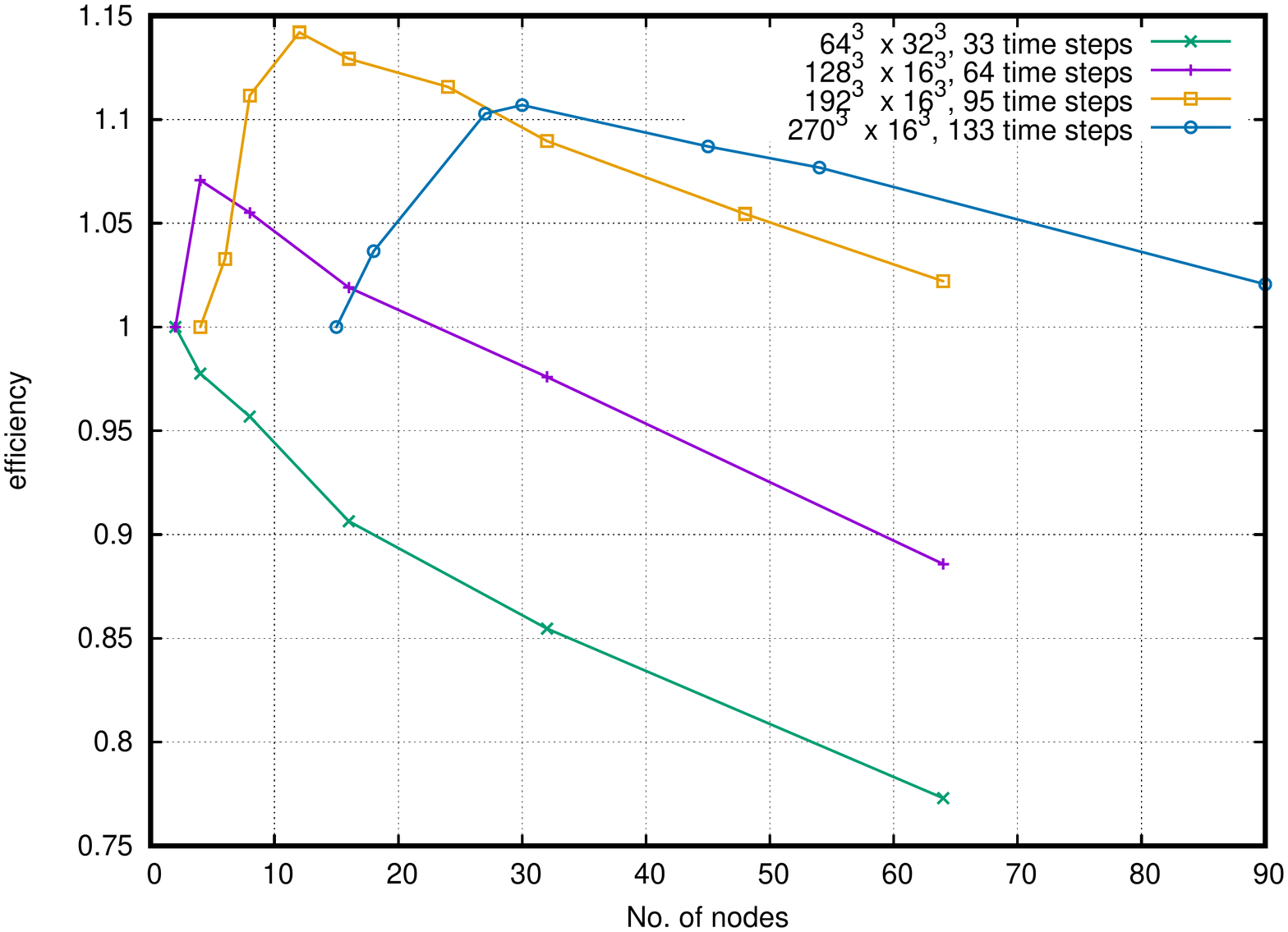}}}
  }
  \caption{Speedup (left) and efficiency (right) as a function of number of computational nodes for BGK. A
    minimal number of nodes requisite to run the test
    is equal to 2 for $64^3 \times 32^3$ and $128^3
    \times 16^3$ mesh and to 4 for the $192^3 \times 16^3$ mesh.}
  \label{fig:speedup}
\end{figure}

Let us first discuss numerical results for the BGK collision kernel.
Simulations were run on four different meshes:
\begin{itemize}
\item $N= 64^3$  and $N_v = 32^3$,
\item $N= 128^3$ and $N_v = 16^3$,
\item $N= 192^3$ and $N_v = 16^3$.
\item $N= 270^3$ and $N_v = 16^3$.  
\end{itemize}
The memory required for storage of the mass array was equal to $64$Gb
in the first two cases, to $216$Gb in the third case and to $600$Gb in
the last case. The minimal number of computational nodes required to
run the Sod test case was therefore equal to $2$ for the first two
cases, to $4$ in the third case and to $15$ in the last case.

Let us first compare different domain decomposition strategies for the
first mesh ($N=62^3$ and $N_v = 32^3$) to find the optimal one. The
computations were run on $N_{\MPI}=2$, $4$, $8$, $16$, $32$ and $64$
nodes. The domain was decomposed into slabs, pencils or cuboids. The
elapsed run time is shown in Table \ref{BGK_domain_decomp}. For
$N_{\MPI}=4$ nodes the best results are obtained for the computational
domain divided into 4 slabs. This coincides with the smallest number
of the ghost cells. For 8 nodes the shortest run time is still obtained
for the slab configuration, even though the number of ghost cells is
$28\%$ bigger than for the pencil configuration. The run time for the
latter is only slightly bigger. The cuboid type domain decomposition
gives the biggest run time. For 16 nodes the run time for slab and
pencil like decomposition is very close to each other with relative
difference of the order of $1\%$. We note that the best performance
was obtained for a pencil like decomposition with $N_{\MPI_x}=8$ and
$N_{\MPI_y}=2$. For 32 and 64 nodes the best performance was obtained
again for pencil like decomposition. The cuboid type decomposition
resulted in a run time comparable with the slab configuration and from
$7\%$ to $12\%$ bigger than the run time for the best configuration
despite the fact that the number of ghost cells was the smallest
(almost $5$ time smaller than the number of ghost cells for slab
configuration on 64 nodes). The weak performance of the cuboid
configuration is due to a greater MPI overhead caused by communication
with significantly bigger number of nodes.

The scalability test were performed with the optimal domain
decomposition strategy. The run time as a function of computational
cores is presented on the Figure \ref{fig:cores_runtime} and the
speedup with parallel efficiency (relative to the smallest number of
nodes employed in the test) on the Figure \ref{fig:speedup}. For the
second mesh the scaling is almost perfect with some super-linear
behaviour when passing from 2 to 4 computational nodes and some loss
of performance when going from 32 to 64 for nodes. For the smallest
mesh there is no super-linearity observed, the performance loss is
observed when passing from 16 to 32 and then to 64 computational
nodes. For the third and fourth mesh the super-linearity appears
between 4 and 16 nodes employed, parallel efficiency is above $0.95$
for 2--32 nodes. The scaling for the last mesh seems to be super
linear in the whole tested range with parallel efficiency above
one. This super linear scaling is probably due to the CPU cache
performance. For all four mesh sizes there is a linear decrease in
parallel efficiency for growing number of nodes.

\begin{table}[t]
  \small
  \centering  
  \begin{tabular}{cccccccc}
    \hline
    \hline
    $N_{\MPI}$ & $N_{\MPI_x}$ & $N_{\MPI_y}$& $N_{\MPI_z}$ &  time       &  \#neigh.&  \#cells &  \#ghosts \\
    \hline
    \hline
    2 & 2 & 1 & 1  &  $\bm{1376.51}$  &     2   &  131072 &  \bf{8192}\\
    \hline                     
    \multirow{2}{*}{4} & 4 & 1 & 1  &  $\bf{704.064}$  &     2   &   65536 &  \bf{8192}\\
    & 2 & 2 & 1  &  $712.356$  &     8   &   65536 &  8448\\
    \hline                     
    \multirow{3}{*}{8} & 8 & 1 & 1  &  $\bf{359.652}$  &     2   &   32768 &  8192\\
    & 4 & 2 & 1  &  $362.894$  &     8   &   32768 &  \bf{6400}\\
    & 2 & 2 & 2  &  $386.392$  &    26   &   32768 &  6536\\
    \hline                     
    \multirow{4}{*}{16}& 16& 1 & 1  &  $191.856$  &     2   &   16384 &  8192\\
    & 8 & 2 & 1  &  $\bf{189.828}$  &     8   &   16384 &  5376\\
    & 4 & 4 & 1  &  $191.113$  &     8   &   16384 &  \bf{4352}\\
    & 4 & 2 & 2  &  $201.823$  &    26   &   16384 &  4424\\
    \hline                     
    \multirow{5}{*}{32}& 32& 1 & 1  &  $106.414$  &     2   &    8192 &  8192\\
    & 16& 2 & 1  &  $101.946$  &     8   &    8192 &  4864\\
    & 8 & 4 & 1  &  $\bf{100.669}$  &     8   &    8192 &  3328\\
    & 8 & 2 & 2  &  $106.942$  &    26   &    8192 &  3368\\
    & 4 & 4 & 2  &  $108.903$  &    26   &    8192 &  \bf{2824}\\
    \hline                     
    \multirow{7}{*}{64}& 64& 1 & 1  &  $63.5180$  &     2   &    4096 &  8192\\
    & 32& 2 & 1  &  $58.7658$  &     8   &    4096 &  4608\\
    & 16& 4 & 1  &  $\bf{55.6507}$  &     8   &    4096 &  2816\\
    & 16& 2 & 2  &  $59.4585$  &    26   &    4096 &  2840\\
    & 8 & 8 & 1  &  $56.5515$  &     8   &    4096 &  2304\\
    & 8 & 4 & 2  &  $59.3814$  &    26   &    4096 &  2024\\
    & 4 & 4 & 4  &  $62.6682$  &    26   &    4096 &  \bf{1736} \\
    \hline 
    \hline
  \end{tabular}
  \caption{Computational time, number of neighboring nodes, number of
    data cells per node and number of ghost cells per node as a
    function of MPI size ($N_{\MPI}$) and domain decomposition (MPI
    dimensions $N_{\MPI_x}$, $N_{\MPI_y}$, $N_{\MPI_z}$) for the BGK kernel test case on
    the $64^3 \times 32^3$ mesh. The smallest number of ghost cells
    and the smallest run time are written in bold type for every tested
    MPI size $N_{\MPI}$.}
    \label{BGK_domain_decomp}
\end{table}

\begin{table}[t]
  \small
  \centering
  \begin{tabular}{ccccccccc}
    \hline
    \hline
    $N_V $ & {\begin{sideways} Vel.\end{sideways}}
    & Cell \# & \#nodes & $N_{\text{cycle}}$ & Time(s) & $T_{\text{cycle}}$ & $T_{\text{cell}}$ & $T_{\text{cell/node}}$
    \\
    &
    \\
    \hline
    \hline
    & 
    \multirow{7}{*}{\begin{sideways}  $[-15,15]$ \end{sideways}}
    && 15 (300 cores)  && $6815.39$ & $51.2$ & $2.60\ 10^{-6}$ & $3.91\ 10^{-5}$    \\
    &&& 18 (360 cores) && $5479.12$ & $41.2$ & $2.09\ 10^{-6}$ & $3.77\ 10^{-5}$    \\
    &&& 27 (540 cores) && $3433.45$ & $25.8$ & $1.31\ 10^{-6}$ & $3.54\ 10^{-5}$    \\
    $16^3$&& $270^3\times 16^3$& 30 (600 cores) &$133$ & $3078.33$ & $23.1$ & $1.18\ 10^{-6}$ & $3.53\ 10^{-5}$ \\
    &&$=80.6\times 10^9$ & 45 (900 cores)  && $2089.94$ & $15.7$ & $7.98\ 10^{-7}$ & $3.59\ 10^{-5}$    \\
    &&& 54 (1080 cores) && $1758$    & $13.2$ & $6.72\ 10^{-7}$ & $3.63\ 10^{-5}$    \\
    &&& 90 (1800 cores) && $1112.93$ & $8.37$ & $4.25\ 10^{-7}$ & $3.83\ 10^{-5}$    \\
    \hline
    \hline            
  \end{tabular}
  \caption{Performance tests on $270^3\times 16^3$ mesh for BGK. Time per
    cycle is obtained by $T_\text{cycle} = T/N_\text{cycle}$, time per
    cycle per cell by $T_\text{cell} = T_\text{cycle}/N_c$ and time
    per cycle per node by $T_\text{cycle/node} =  N_s T_\text{cell}$.}
  \label{tab:270}
\end{table}

\begin{table}[t]
  \small
  \centering  
  \begin{tabular}{ccccccccc}
    \hline
    \hline
    $N_V $ & {\begin{sideways} Vel.\end{sideways}}
    & Cell \# & \#nodes & $N_{\text{cycle}}$ & Time(s) & $T_{\text{cycle}}$ & $T_{\text{cell}}$ & $T_{\text{cycle/node}}$
    \\
    &
    \\
    \hline
    \hline    
    &
    \multirow{9}{*}{\begin{sideways}  $[-15,15]$ \end{sideways}} 
    && 4 (80 cores)    && $6696.82$ & $70.5$ & $9.96\ 10^{-6}$ & $3.98\ 10^{-5}$    \\
    &&& 6 (120 cores)  && $4322.4 $ & $45.5$ & $6.43\ 10^{-6}$ & $3.86\ 10^{-5}$    \\
    &&& 8 (160 cores)  && $3012.17$ & $31.7$ & $4.48\ 10^{-6}$ & $3.58\ 10^{-5}$    \\
    &&& 12 (240 cores) && $1954.71$ & $20.6$ & $2.91\ 10^{-6}$ & $3.49\ 10^{-5}$    \\
    $16^3$&&$192^3\times 16^3$&16 (320 cores) &$95$& $1482.56$ & $15.6$ & $2.20\ 10^{-6}$ & $3.53\ 10^{-5}$ \\
    &&$=29\times 10^9$    & 24 (480 cores)  && $1000.37$ & $10.5$ & $1.49\ 10^{-6}$ & $3.57\ 10^{-5}$    \\
    &&& 32 (640 cores)  && $768.273$ & $8.09$ & $1.14\ 10^{-6}$ & $3.66\ 10^{-5}$    \\
    &&& 48 (960 cores)  && $529.223$ & $5.57$ & $7.87\ 10^{-7}$ & $3.78\ 10^{-5}$    \\
    &&& 64 (1280 cores) && $409.52 $ & $4.31$ & $6.09\ 10^{-7}$ & $3.90\ 10^{-5}$    \\
    \hline
    \hline
  \end{tabular}
  \caption{Performance tests on $192^3\times 16^3$ mesh for BGK. Time per
    cycle is obtained by $T_\text{cycle} = T/N_\text{cycle}$, time per
    cycle per cell by $T_\text{cell} = T_\text{cycle}/N_c$ and time
    per cycle per node by $T_\text{cycle/node} =  N_s T_\text{cell}$.}
  \label{tab:192}
\end{table}

\begin{table}[!h]
  \small
  \centering  
  \begin{tabular}{ccccccccc}
    \hline
    \hline
    $N_V $ & {\begin{sideways} Vel.\end{sideways}}
    & Cell \# & \#nodes & $N_{\text{cycle}}$ & Time(s) & $T_{\text{cycle}}$ & $T_{\text{cell}}$ & $T_{\text{cell/node}}$
    \\
    &
    \\    
    \hline
    \hline
    &
    \multirow{7}{*}{\begin{sideways}  $[-15,15]$ \end{sideways}}     
    && 2 (40 cores)   && $2499.49$ & $39.1$ & $1.86\ 10^{-5}$ & $3.72\ 10^{-5}$    \\
    &&& 4 (80 cores)  && $1167.2 $ & $18.2$ & $8.70\ 10^{-6}$ & $3.48\ 10^{-5}$    \\
    \multirow{2}{*}{$16^3$}&&$128^3\times 16^3$& 8 (160 cores) & \multirow{2}{*}{$64$} & $592.256$ & $9.25$ & $4.41\ 10^{-6}$ & $3.53\ 10^{-5}$\\
    &&$=8.6\times 10^9$&16 (320 cores)  && $306.597$ & $4.79$ & $2.28\ 10^{-6}$ & $3.65\ 10^{-5}$    \\
    &&& 32 (640 cores)     && $160.077$ & $2.50$ & $1.19\ 10^{-6}$ & $3.82\ 10^{-5}$    \\
    &&& 64 (1280 cores)    && $88.1792$ & $1.38$ & $6.57\ 10^{-7}$ & $4.20\ 10^{-5}$    \\
    \hline
    \hline
  \end{tabular}
  \caption{Performance tests on $128^3\times 16^3$ mesh for BGK. Time per
    cycle is obtained by $T_\text{cycle} = T/N_\text{cycle}$, time per
    cycle per cell by $T_\text{cell} = T_\text{cycle}/N_c$ and time
    per cycle per node by $T_\text{cycle/node} =  N_s T_\text{cell}$.}
  \label{tab:128}
\end{table}

Tables \ref{tab:270}--\ref{tab:64} show the total run time $T$, time
spent on a single cycle $T_\text{cycle} = T/N_\text{cycle}$, average
time time spent on a single cell per cycle $T_\text{cell} =
T_\text{cycle}/N_c$ and finally the average time spent by a single
node on a single cell per cycle $T_\text{cell/node} = N_{\MPI}
T_\text{cycle}$ for different number of nodes employed. The
performance loss is observed when the subdomain size become relatively
small. This is manifested in the increase of the average time per
cycle spent by one node on a single cell.  This is due to the fact
that the internodal communication time becomes comparable with the
time spent on computations.

The size of spatial mesh allocated in every node seems also to
influence the cache performance at the nodes. On Tables \ref{tab:270},
\ref{tab:192} and \ref{tab:128} the non linear scaling is observed for
small number of nodes --- the average time per cell per cycle per node
is decreasing when the number of nodes increases. The best performance
is observed on 30 nodes for the $270^3$ mesh ($656\ 10^3$ cells per
node), on 12 nodes for $192^3$ mesh ($590\ 10^3$ cells per node) and
on 4 nodes on $128 ^3$ mesh ($524\ 10^3$ cells per node).

There is no such effect observed on Table \ref{tab:64}. Let us now
compare Tables \ref{tab:128} and \ref{tab:64}. The dimension of the
problem is the same. In the first case the domain was discretized with
$128^3$ points in the physical space and with the $16^3$ velocity
points. In the second case a spatial mesh of the size $64^3$ was used
and a velocity mesh of the size $32^3$, giving a total of $8.6\ 10^9$
degrees of freedom. The size of the spatial mesh per node is eight
times smaller in second case while the number of degrees of freedom is
the same in both cases. The $N=128^3$ mesh seems to perform better
than the $64^3$ one with time per cycle being smaller and with better
MPI efficiency. This is not surprising as the amount of data exchanged
with neighbouring nodes is at least 2 times smaller.

\begin{table}[t]
  \small
  \centering  
  \begin{tabular}{ccccccccc}
    \hline
    \hline
    $N_V $ & {\begin{sideways} Vel.\end{sideways}}
    & Cell \# & \#nodes & $N_{\text{cycle}}$ & Time(s) & $T_{\text{cycle}}$ & $T_{\text{cell}}$ & $T_{\text{cell/node}}$
    \\
    &    
    \\    
    \hline
    \hline
    & 
    \multirow{7}{*}{\begin{sideways}  $[-15,15]$ \end{sideways}}
    && 2 (40 cores)     && $1250.51$ & $37.9$ & $1.45\ 10^{-4}$ & $2.89\ 10^{-4}$    \\
    &&& 4 (80 cores)    && $635.297$ & $19.3$ & $7.34\ 10^{-5}$ & $2.94\ 10^{-4}$    \\
    \multirow{2}{*}{$32^3$}&&$64^3\times 32^3$& 8 (160 cores)   &\multirow{2}{*}{$33$}& $328.987$ & $9.97$ & $3.80\ 10^{-5}$ & $3.04\ 10^{-4}$ \\
    &&$=8.6\times 10^9$& 16 (320 cores)   && $172.168$ & $5.22$ & $1.99\ 10^{-5}$ & $3.18\ 10^{-4}$    \\
    &&& 32 (640 cores) && $98.2927$ & $2.98$ & $1.14\ 10^{-5}$ & $3.64\ 10^{-4}$    \\
    &&& 64 (1280 cores) && $58.6773$ & $1.78$ & $6.78\ 10^{-6}$ & $4.34\ 10^{-4}$   \\
    \hline
    \hline
  \end{tabular}
  \caption{Performance tests on $64^3\times 32^3$ mesh for BGK. Time per
    cycle is obtained by $T_\text{cycle} = T/N_\text{cycle}$, time per
    cycle per cell by $T_\text{cell} = T_\text{cycle}/N_c$ and time
    per cycle per node by $T_\text{cycle/node} =  N_s T_\text{cell}$.}
  \label{tab:64}
\end{table}
\begin{table}[h]
  \small
  \centering    
  \begin{tabular}{|cc|c|c||c||cc|}
    \hline
    &&  \multirow{2}{*}{\textbf{Cycle}} &  \textbf{CPU}  &  \multirow{2}{*}{\textbf{Main routines}}
    &  \textbf{Cost CPU }  &  \textbf{Cost}  \\
    && & (s) & & vs total (s) & vs total (\%)  \\
    \hline
    \hline
    \multirow{6}{*}{\begin{sideways} \textbf{$64^3 \times 32^3$}  \end{sideways}}
    &
    \multirow{6}{*}{\begin{sideways} \textbf{2 nodes}  \end{sideways}}
    & \multirow{6}{*}{ 32 }
    & \multirow{6}{*}{ 1374 }
    &Transport          &   0.0028 & $2\ 10^{-4}$\%  \\
    \cline{5-7}
    && & &ToConservative &  269.123 & 19.59\%  \\
    \cline{5-7}
    && & &ToPrimitive    &  0.01  & $7\ 10^{-4}$\%  \\
    \cline{5-7}
    && & &Collision     &  904.166 & 65.8\% \\  
    \cline{5-7}
    && & &Communication &  200.718 & 14.61\% \\  
    \cline{5-7}    
    \cline{5-7}
    && & & = & 1374.02 & 100\%  \\
    \hline
    \hline
    \multirow{6}{*}{\begin{sideways} \textbf{$64^3 \times 32^3$}  \end{sideways}}
    &
    \multirow{6}{*}{\begin{sideways} \textbf{64 nodes}  \end{sideways}}
    & \multirow{6}{*}{ 32 }
    & \multirow{6}{*}{ 56 }
    &Transport          &  0.0027 & $4.9\ 10^{-3}$\%  \\
    \cline{5-7}
    && & &ToConservative &  18.06 &  32.34\%  \\
    \cline{5-7}
    && & &ToPrimitive    &  0.0008 &  $1.4\ 10^{-3}$\%  \\
    \cline{5-7}
    && & &Collision     & 28.113 & 50.34\%  \\  
    \cline{5-7}
    && & &Communication &  9.67 & 17.31\% \\  
    \cline{5-7}    
    \cline{5-7}
    && & & = & 55.85 & 100\%  \\
    \hline
    \hline
    \multirow{6}{*}{\begin{sideways} \textbf{$270^3 \times 16^3$}  \end{sideways}}
    &
    \multirow{6}{*}{\begin{sideways} \textbf{90 nodes}  \end{sideways}}
    & \multirow{6}{*}{ 133 }
    & \multirow{6}{*}{ 1038 }
    &Transport          &  0.0041 & $4\ 10^{-4}$\%  \\
    \cline{5-7}
    && & &ToConservative &  234.78 &  22.62\%  \\
    \cline{5-7}
    && & &ToPrimitive    &  0.053 &  $5\ 10^{-3}$\%  \\
    \cline{5-7}
    && & &Collision     & 751.36 & 72.4\%  \\  
    \cline{5-7}
    && & &Communication &  51.59 & 4.97\% \\  
    \cline{5-7}    
    \cline{5-7}
    && & & = & 1037.76 & 100\%  \\
    \hline    
  \end{tabular}  
  \caption{ \label{tab:profiling3D-BGK} Profiling of the average cost
    for each routine for the BGK collision operator on $64^3 \times
    32^3$ mesh for 32 time steps and for $270^3 \times 16^3$ mesh for
    133 time steps.  }
\end{table}

Table \ref{tab:profiling3D-BGK} shows profiling data for the $64^3
\times 32^3$ mesh on 2 and 64 nodes as well as for the $270^3\times
16^3$ mesh on 90 computational nodes. For the coarse mesh the time
spent on the BGK collision operator decreases from $65\%$ (on 2 nodes)
to $50\%$ (on 64 nodes) while on the fine mesh the relaxation term
takes more than $70\%$ of the computational time (on 90 nodes). This
explains the efficiency loss for larger number of nodes: the time
spent on the BGK operator is too small to hide internodal
communication. Please note that better parallelization strategy can be
applied for this collision model: decomposition of the velocity space
(instead of physical) would require less data exchange between MPI
processes in this particular case and yield therefore a better
parallel efficiency. This alternative approach is however not suited
for the Boltzmann collision operator.

\subsection{3D Boltzmann}

Let us now turn our attention to the 3D Boltzmann collision
kernel. This task is much more demanding in terms of computational
time as the relaxation routine involves multiple and expensive calls
to the Fast Fourier Transform. The tests were run on following meshes:
\begin{itemize}
\item $N= 64^3$ and $N_v = 16^3$,
\item $N= 64^3$ and $N_v = 32^3$.
\item $N= 128^3$ and $N_v = 32^3$. 
\end{itemize}
The tests for the $N=64^3$ meshes were performed on slab type domain
decomposition.
The results are presented on Figures \ref{fig:cores_runtimeBoltz} and
\ref{fig:speedupBoltz} and on Tables \ref{tab:64x16B} and
\ref{tab:64x32B}. For the last mesh computations were performed for
$16$ time steps and until the final time $t_f = 0.07$ was
reached. This clearly does not allow to present any physically
interesting results but is enough to study the parallel efficiency of
the herein proposed method as during this time every particle will
change a space cell at least once. The tests were run on $N_{\MPI}$
ranging from $16$ to $1024$. That is to say on $128$ to $2048$
processors and on $1024$ to $16384$ computational cores. The results
are presented on Table \ref{tab:128x32B} and on Figures
\ref{fig:cores_runtimeBoltz} and \ref{fig:speedupBoltz}.

A comparison of the run times obtained for the second mesh (see Figure
\ref{fig:cores_runtimeBoltz}) with the run time obtained for the BGK
collision kernel shows that the latter is approximately 100 times
faster. As the time spent on the computation of the relaxation kernel
is now much more important, the communication time becomes negligible
even if the size of the subdomains is small. The method shows the
strong scaling that is close to ideal for all tested meshes with no
efficiency loss observed (see Figure \ref{fig:speedupBoltz}) and the
parallel efficiency close to one. Also, the average time spent by one
computational node on one cell per cycle does not depend on number of
nodes employed in the computations. The method clearly enjoys the
strong scaling close to ideal in the tested range.

Comparison of the run time for both meshes shows that multiplying the
number of Fourier modes by 2 in each direction results in a run time
multiplied by 16 instead of 8.  This loss of performance is related to
the computational complexity of the fast spectral solver for the
Boltzmann collision kernel which is of $O(N_v \log N_v)$ as well as to
the cache performance. Each CPU performs 10 (8 on Curie) FFTs in
parallel at the same time. The size of one transform required for the
convolution computation is $32^3$ for the $16^3$ Fourier modes. This
means that the memory required to store 10 vectors containing $32^3$
complex values in double precision format for 10 in-place transforms
is $5$MB, which is less than the cache capacity. For the finer mesh
the memory required is $40$MB, that is to say four times more than the
cache of the CPU.

\def\xxxa{0.495\textwidth} 
\begin{figure}[t] 
  \centering
  \includegraphics[width=\xxxa]{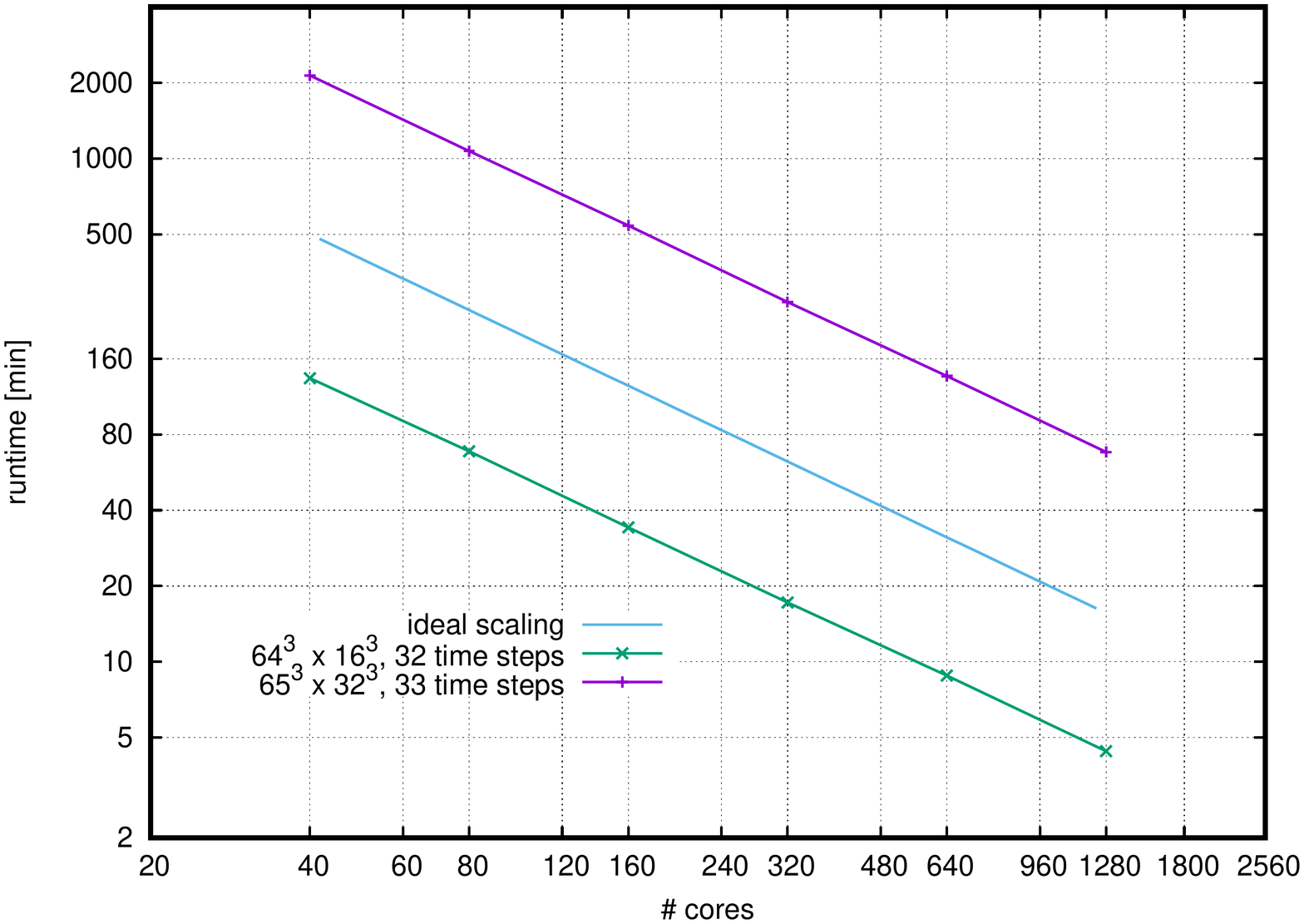}
  \includegraphics[width=\xxxa]{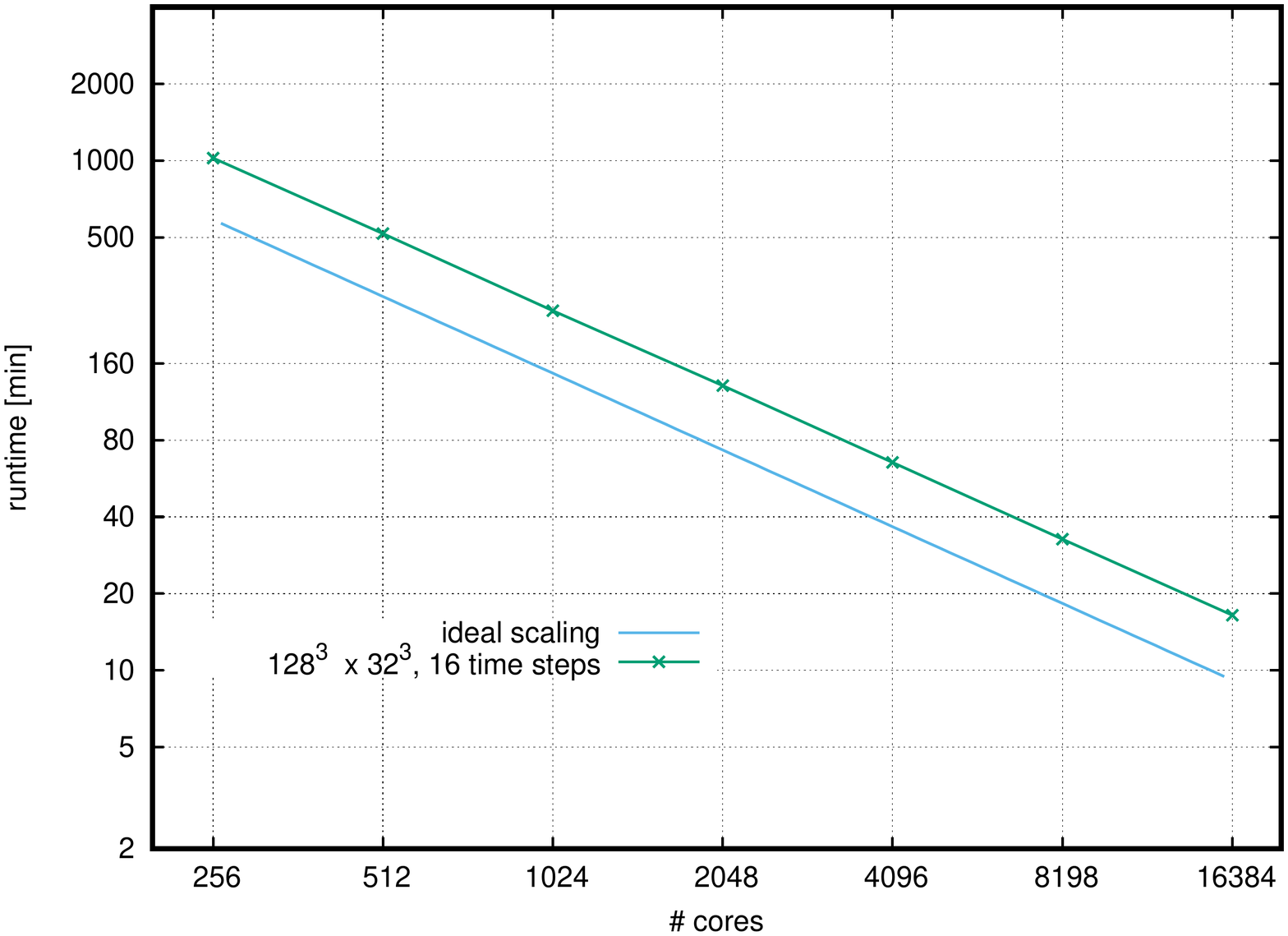}
  \caption{Computational time as a function of number of cores employed on EOS machine (left panel) and on TGCC-CURIE (right panel).}
  \label{fig:cores_runtimeBoltz}
\end{figure}

\begin{figure}[h] 
  \centering
  \includegraphics[width=\xxxa]{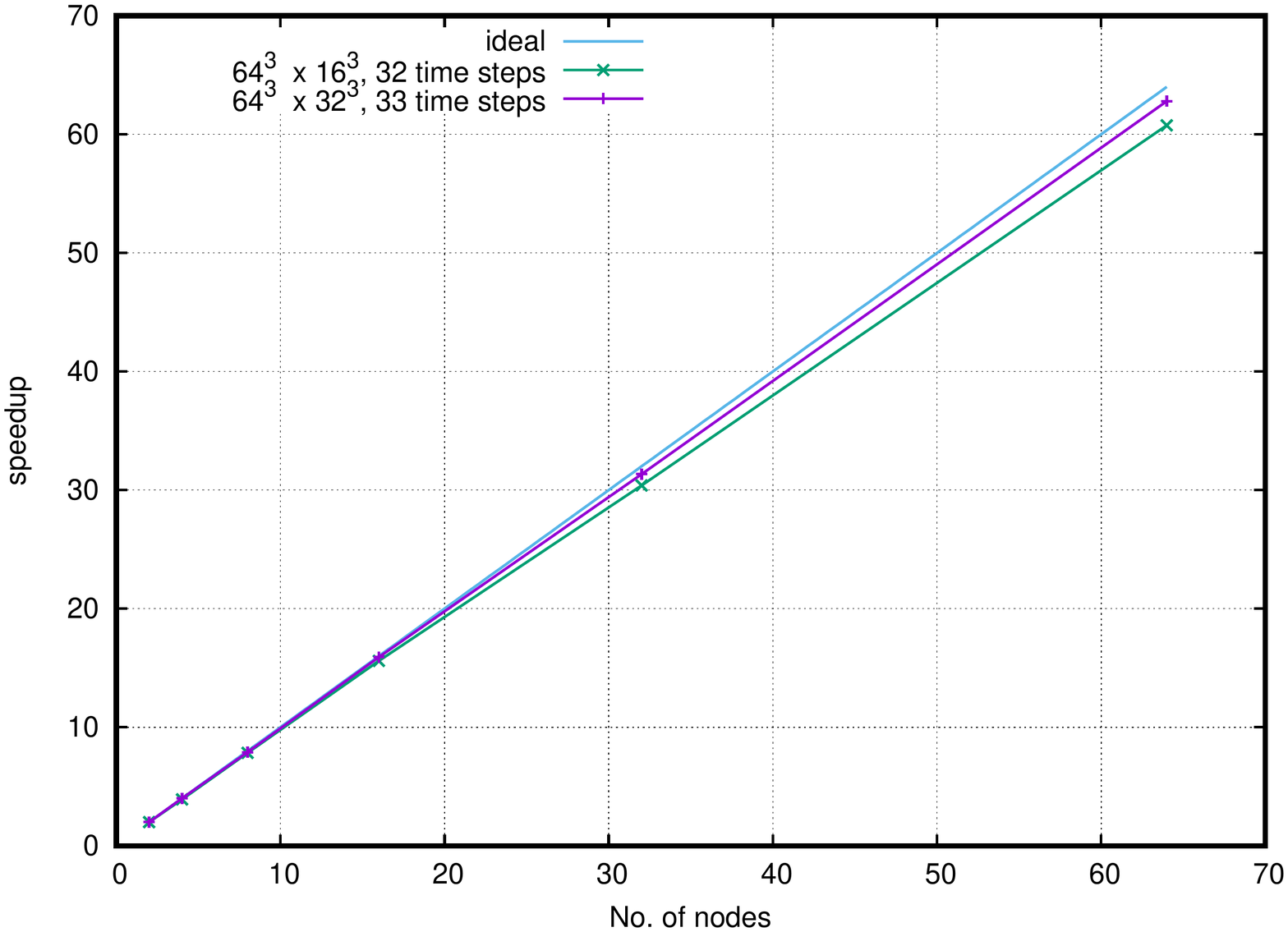}
  \includegraphics[width=\xxxa]{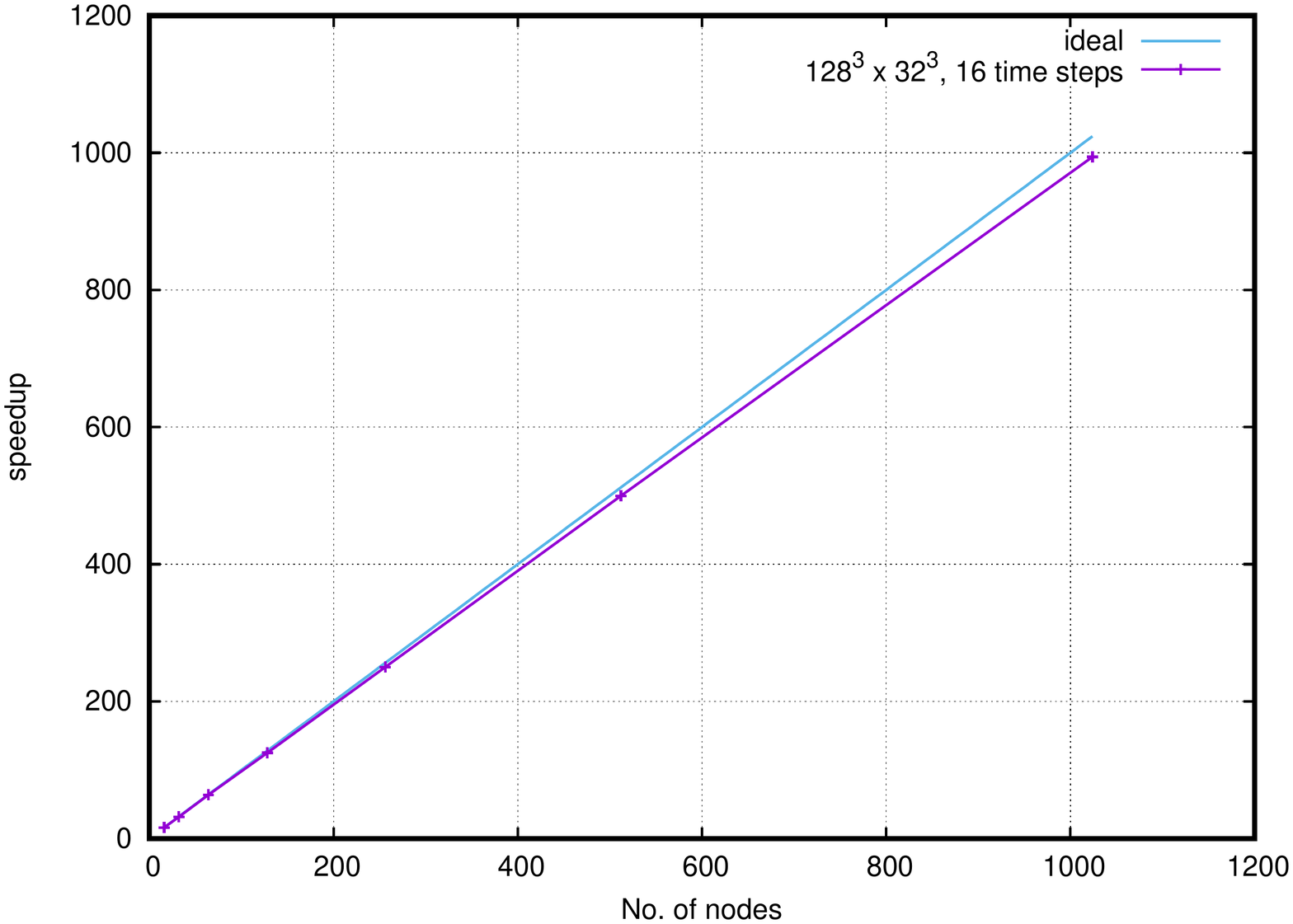}\\
  \includegraphics[width=\xxxa]{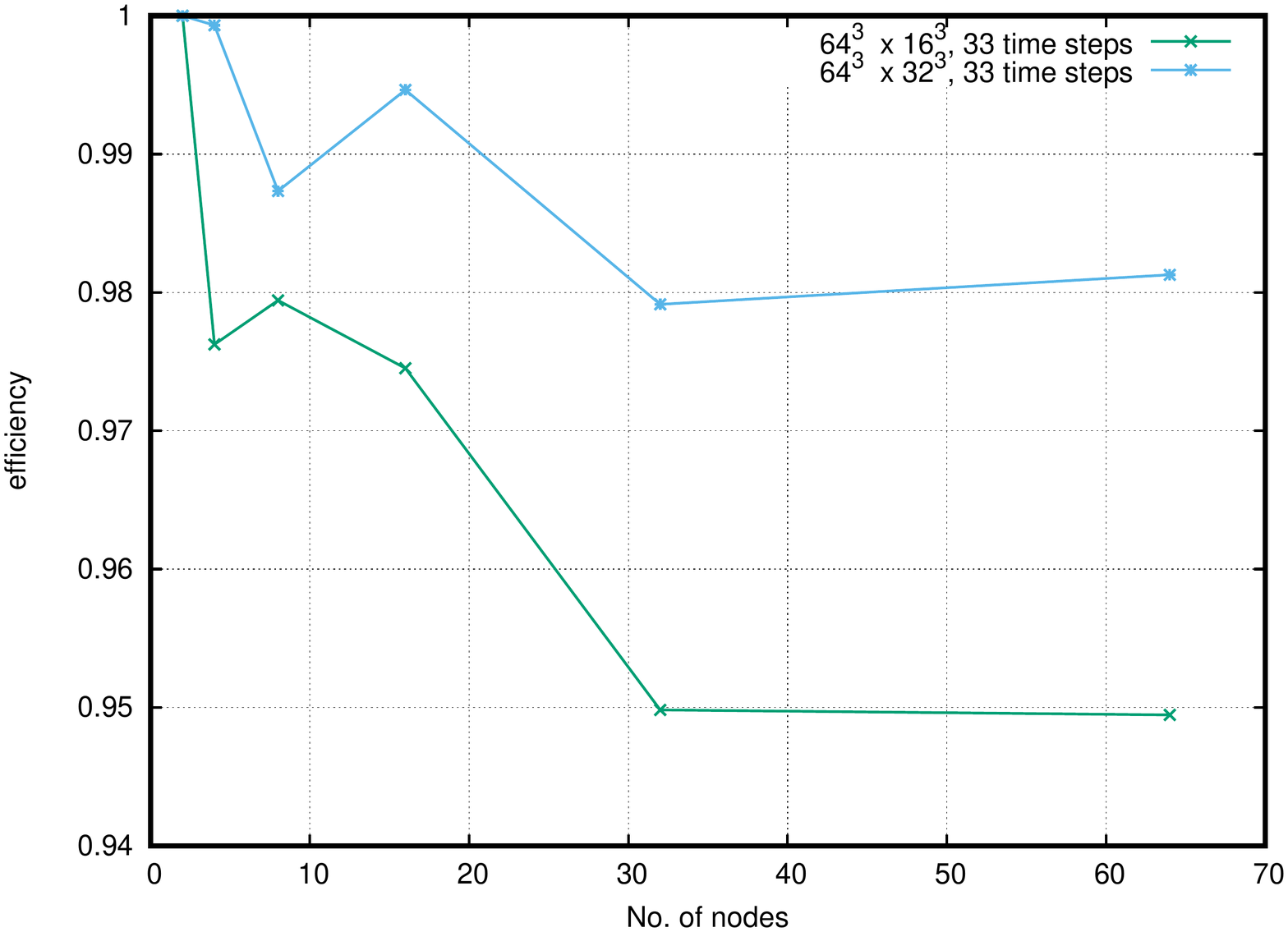}
  \includegraphics[width=\xxxa]{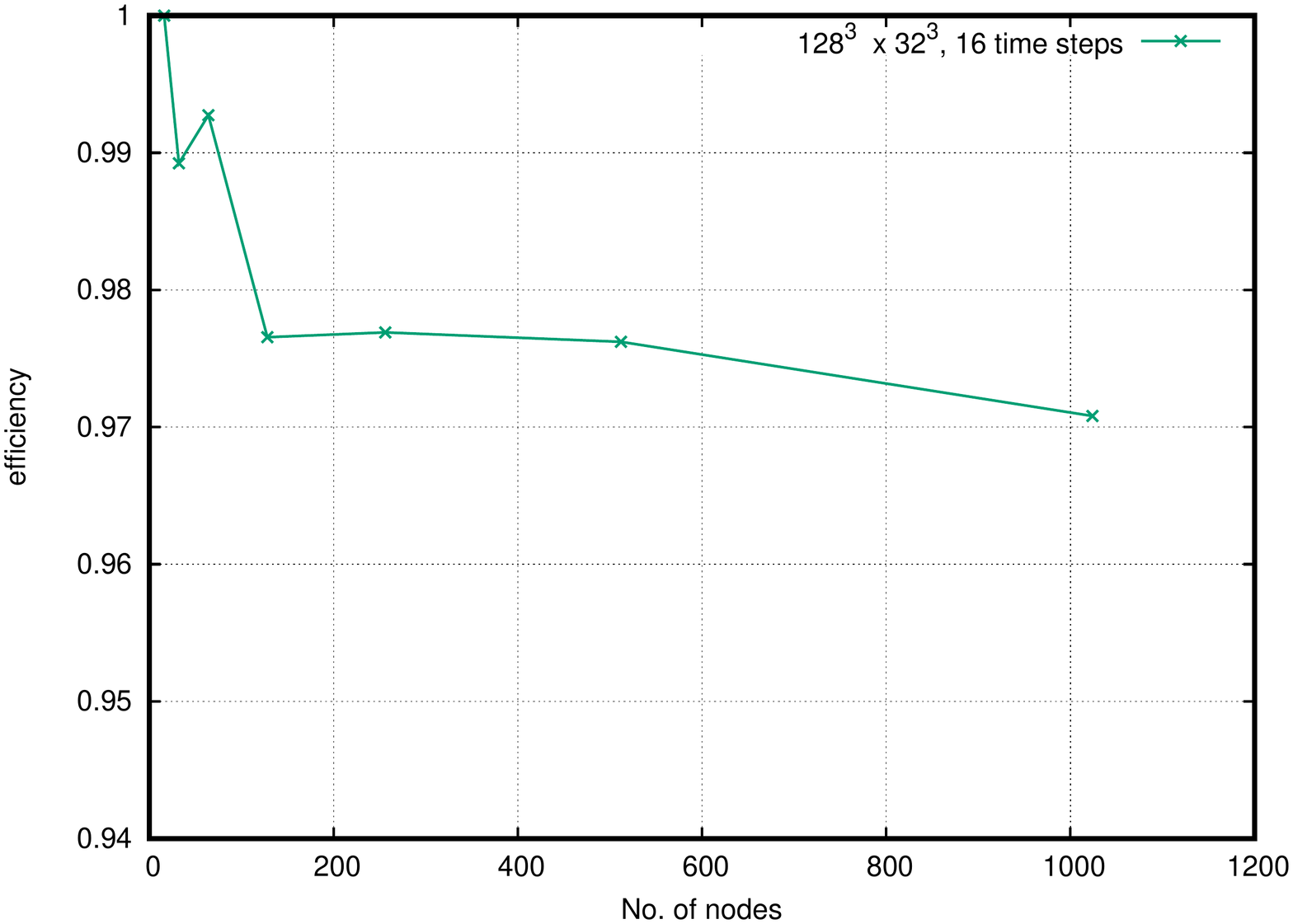}  
  \caption{Speedup (top) and efficienvy (bottom) as a function of number of computational nodes for
    3D Boltzmann collision kernel on EOS machine (left panel) and on TGCC-CURIE (right panel).}
  \label{fig:speedupBoltz}
\end{figure}

\begin{table}[h]
  \small
  \centering  
  \begin{tabular}{ccccccccc}
    \hline
    \hline
    $N_V $ & {\begin{sideways} Vel.\end{sideways}}
    & Cell \# & \#nodes & $N_{\text{cycle}}$ & Time(s) & $T_{\text{cycle}}$ & $T_{\text{cell}}$ & $T_{\text{cell/node}}$
    \\
    &
    \\    
    \hline
    \hline
    & 
    \multirow{7}{*}{\begin{sideways}  $[-15,15]$ \end{sideways}}
    && 2 (40 cores)     && $8043.07$ & $244$ & $9.3\ 10^{-4}$ & $1.86\ 10^{-3}$    \\
    &&& 4 (80 cores)    && $4119.41$ & $125$ & $4.76\ 10^{-4}$ & $1.90\ 10^{-3}$    \\
    \multirow{2}{*}{$16^3$}&&$64^3\times 16^3$& 8 (160 cores)   &\multirow{2}{*}{$33$}& $2053.01$ & $62.2$ & $2.37\ 10^{-4}$ & $1.90\ 10^{-3}$ \\
    &&$=1.07\times 10^9$& 16 (320 cores) && $1031.67$ & $31.3$ & $1.19\ 10^{-4}$ & $1.91\ 10^{-3}$  \\
    &&& 32 (640 cores)  && $529.246$ & $16$ & $6.12\ 10^{-5}$ & $1.96\ 10^{-3}$    \\
    &&& 64 (1280 cores) && $264.725$ & $8.02$ & $3.06\ 10^{-5}$ & $1.96\ 10^{-3}$    \\
    \hline
    \hline
  \end{tabular}
  \caption{Performance tests on $64^3\times 16^3$ mesh for Boltzmann
    collision kernel. Time per cycle is obtained by $T_\text{cycle} =
    T/N_\text{cycle}$, time per cycle per cell by $T_\text{cell} =
    T_\text{cycle}/N_c$ and time per cycle per node by
    $T_\text{cycle/node} = N_s T_\text{cell}$.}
  \label{tab:64x16B}
\end{table}

\begin{table}[h]
  \small
  \centering  
  \begin{tabular}{ccccccccc}
    \hline
    \hline
    $N_V $ & {\begin{sideways} Vel.\end{sideways}}
    & Cell \# & \#nodes & $N_{\text{cycle}}$ & Time(s) & $T_{\text{cycle}}$ & $T_{\text{cell}}$ & $T_{\text{cell/node}}$
    \\
    &    
    \\    
    \hline
    \hline
    & 
    \multirow{7}{*}{\begin{sideways}  $[-15,15]$ \end{sideways}}
    && 2 (40 cores)     && $128378$ & $3890$ & $1.48\ 10^{-2}$ & $2.97\ 10^{-2}$    \\
    &&& 4 (80 cores)    && $64233.2$ & $1950$ & $7.43\ 10^{-3}$ & $2.97\ 10^{-2}$    \\
    \multirow{2}{*}{$32^3$}&&$64^3\times 32^3$& 8 (160 cores)&\multirow{2}{*}{$33$}& $32506.4$ & $985$ & $3.76\ 10^{-3}$ & $3.01\ 10^{-2}$ \\
    &&$=8.6\times 10^9$& 16 (320 cores) && $16133.7$ & $489$ & $1.87\ 10^{-3}$ & $2.98\ 10^{-2}$    \\
    &&& 32 (640 cores)    && $8194.58$ & $248$ & $9.47\ 10^{-4}$ & $3.03\ 10^{-2}$    \\
    &&& 64 (1280 cores)    && $4088.34$ & $124$ & $4.73\ 10^{-4}$ & $3.02\ 10^{-2}$    \\
    \hline
    \hline
  \end{tabular}
  \caption{Performance tests on $64^3\times 32^3$ mesh for Boltzmann
    collision kernel. Time per cycle is obtained by $T_\text{cycle} =
    T/N_\text{cycle}$, time per cycle per cell by $T_\text{cell} =
    T_\text{cycle}/N_c$ and time per cycle per node by
    $T_\text{cycle/node} = N_s T_\text{cell}$.}
  \label{tab:64x32B}
\end{table}

\begin{table}[h]
  \small
  \centering  
  \begin{tabular}{ccccccccc}
    \hline
    \hline
    $N_V $ & {\begin{sideways} Vel.\end{sideways}}
    & Cell \# & \#nodes & $N_{\text{cycle}}$ & Time(s) & $T_{\text{cycle}}$ & $T_{\text{cell}}$ & $T_{\text{cell/node}}$
    \\
    &    
    \\    
    \hline
    \hline
    & 
    \multirow{8}{*}{\begin{sideways}  $[-15,15]$ \end{sideways}}
    && 16 (256 cores)     && $61346.4$ & $3834.2$ & $1.83\ 10^{-3}$ & $2.93\ 10^{-2}$    \\
    &&& 32 (512 cores)    && $31006.9$ & $1937.9$ & $9.24\ 10^{-4}$ & $2.96\ 10^{-2}$    \\
    &&& 64 (1024 cores)   && $15448.8$ & $965.6$ & $4.60\ 10^{-4}$ & $2.95\ 10^{-2}$    \\
    $32^3$&&$128^3\times 32^3$& 128 (2048 cores) &$16$ & $7852.4$ & $490.8$ & $2.34\ 10^{-4}$ & $3.00\ 10^{-2}$    \\
    &&$=69\times 10^{9}$& 256 (4096 cores) && $3924.8$ & $245.3$ & $1.17\ 10^{-4}$ & $2.99\ 10^{-2}$    \\
    &&& 512 (8192 cores)   && $1963.7$ & $122.7$ & $5.85\ 10^{-5}$ & $3.00\ 10^{-2}$    \\
    &&& 1024 (16384 cores) && $987.37$ & $61.71$ & $2.94\ 10^{-5}$ & $3.01\ 10^{-2}$    \\
    \hline
    \hline
  \end{tabular}
  \caption{Performance tests on $128^3\times 32^3$ mesh for Boltzmann
    collision kernel for 16 time steps. Time per cycle is obtained by $T_\text{cycle} =
    T/N_\text{cycle}$, time per cycle per cell by $T_\text{cell} =
    T_\text{cycle}/N_c$ and time per cycle per node by
    $T_\text{cycle/node} = N_s T_\text{cell}$.}
  \label{tab:128x32B}
\end{table}


\begin{table}
  \small
  \centering  
  \begin{tabular}{|c|c|c||c||cc|}
    \hline
    &  \multirow{2}{*}{\textbf{Cycle}} &  \textbf{CPU}  &  \multirow{2}{*}{\textbf{Main routines}}
    &  \textbf{Cost CPU }  &  \textbf{Cost}  \\
    & & (s) & & vs total (s) & vs total (\%)  \\
    \hline
    \hline
     \multirow{6}{*}{\begin{sideways} \textbf{16 nodes}  \end{sideways}}
    & \multirow{6}{*}{ 16 }
    & \multirow{6}{*}{ 61346 }
    &Transport          &   0.0032 & $5\ 10^{-6}$\%  \\
    \cline{4-6}
    & & &ToConservative &  175.2 & 0.29\%  \\
    \cline{4-6}
    & & &ToPrimitive    &  0.0056  & $9\ 10^{-6}$\%  \\
    \cline{4-6}
    & & &Collision     &  60418.6 & 98.5\% \\  
    \cline{4-6}
    & & &Communication &  752.59 & 1.23\% \\  
    \cline{4-6}    
    \cline{4-6}
    & & & = & 61346.4 & 100\%  \\
    \hline
    \hline
    \multirow{5}{*}{\begin{sideways} \textbf{1024 nodes}  \end{sideways}}
    & \multirow{5}{*}{ 16 }
    & \multirow{5}{*}{ 987 }
    &Transport          &  0.0043 & $4.3\ 10^{-3}$\%  \\
    \cline{4-6}
    & & &ToConservative &  6.16 &  0.62\%  \\
    \cline{4-6}
    & & &ToPrimitive    &  0.00028 &  $2.9\ 10^{-6}$\%  \\
    \cline{4-6}
    & & &Collision     & 944.222 & 95.63\%  \\  
    \cline{4-6}
    & & &Communication &  36.98 & 3.75\% \\  
    \cline{4-6}    
    \cline{4-6}
    & & & = & 987.365 & 100\%  \\
    \hline
  \end{tabular}
  \caption{ \label{tab:profiling3Dcurie} Profiling of the average cost
  for each routine for the Boltzmann collision kernel on $128^3 \times
  32^3$ mesh for 16 time steps.
  }
\end{table}

Table \ref{tab:profiling3Dcurie} shows the profiling data for 16 and
1024 computational nodes. Time spent on MPI communication increases
from $1.23\%$ to $3.75\%$. Collision step takes above $95\%$ of the
total run time and hides the parallelization costs.

%
\section{Conclusions}

In this paper we have presented an extension of the Fast Kinetic
Scheme introduced in \cite{FKS} to the Boltzmann collision operator by
means of Fast Spectral Method and possible MPI parallelization. The
obtained strong scaling is closed to ideal in the tested range (up to
1024 computational nodes) for 3D$\times$3D problems for the Boltzmann
operator. In the proposed method the physical space was distributed
over computational nodes by means of the MPI. Each node was supplied
with a complete velocity grid. This approach has proven to be well
suited for velocity non local collision operators (like Boltzmann) and
performs less good for simplified and velocity local models (BGK). The
reason is that the BGK collision kernel is much less resource
demanding --- computational time spent on particle interaction in one
space cell is two orders of magnitude smaller compared to the full
Boltzmann operator. Better performance for the BGK operator could be
obtained by employing the following alternative: distribute velocity
space over computational nodes and supply every node with a complete
physical space. This approach would minimize the communication for the
BGK kernel. However, the optimal algorithm for this operator is not
the scope of this paper and numerical results were presented only to
illustrate the computational complexity of the Boltzmann operator.
The algorithm was tested on classical architectures with collision
kernels parallelized locally by the means of OpenMP. We expect to
maintain the similar scaling on GPU or Intel MIC based node
architectures, which is the scope of future works.


%
%
\section*{Acknowledgements}

This work has been supported by the french 'Agence Nationale pour la
Recherche (ANR)' in the frame of the contract ``Moonrise''
(ANR-11-MONU-009-01). Numerical simulations were performed using HPC
resources from GENCI-TGCC (Grant 2016-AP010610045) and from CALMIP
(project P1542).  The author would like to acknowledge fruitful
discussions with G.Dimarco, F.Filbet and R.Loub\`ere.

%
%
%

\appendix

\section{Boltzmann collision operator}\label{Boltz}

The Boltzmann collision operator is solved by the means of Fast
Spectral Scheme presented in this section for a selected time step
$t^n$ and a selected grid point $x_j$. The same computation is
repeated for other grid points and time steps. In this section we
denote by $f$ the distribution function (of $v$ only) at time $t^n$
and point $x_j$: $f=f(v)=f(x_j,v,t^n)$.

In order to compute the Boltzmann integral (\ref{eq:boltzmann}), let
us suppose that the distribution function $f$ has a compact support on
the ball $B_0(R)$ of radius $R$ centered in the origin. It can be
shown \cite{pareschi:1996,PR2} that the support of the collision operator
$Q(f,f)$ is included in the ball $B_0(\sqrt{2}R)$ and
\begin{gather}
  \Q_B (f,f) = \Int_{B_0(2R)} \Int_{S^2}
  B(|g|,\theta )
  \left(
  f(v^\prime)f(v_1^\prime) - f(v)f(v-g)
  \right)
  d\omega dg
\nonumber 
\end{gather}
with $v^\prime,v_1^\prime,v-g \in B_0((2+\sqrt{2})R)$. We can
therefore restrict $f$ to the cube $[-T,T]^3$ with $T\geq
(2+\sqrt{2})R$ assuming $f(v) =0$ on $[-T,T]^3\setminus B_0(R)$ and
then extend it to a periodic function on $[-T,T]^3$. As a consequence
of the periodicity of $f$ it is sufficient to take $T \geq
(3+\sqrt{2})R/2$ to prevent overlapping of the regions where $f$ is
different from zero \cite{PR2}. In order to simplify notation let us
take $T = \pi$ and $R = \lambda \pi$ with $\lambda
=2/(3+\sqrt{2})$. Let $Q^R_B(f)$ denote the Boltzmann operator with
cut-off. Let us perform a discrete Fourier transform on $f$ obtaining
\begin{gather}
  f_N(v) = \sum_{k=-N/2}^{N/2} \f_k e^{i k \cdot v}
\nonumber 
,
  \\
  \f_k = \frac{1}{(2\pi)^{3}}\int_{[-\pi,\pi]^{3}} f(v)
  e^{-i k \cdot v }\,dv.
\nonumber 
\end{gather}
with $k$ being a multi-index $k=(k_x,k_y,k_z)$ and $N =
(\sqrt[3]{N_v},\sqrt[3]{N_v},\sqrt[3]{N_v})$ is vector containing the
number of velocity discretization points (Fourier modes) in each
direction.  Let us impose that the residue of the collision step is
orthogonal to any trigonometric polynomial of degree less or equal
than $N$ in order to obtain a set of ODEs for coefficients $\f_k$:
\begin{gather}
  \int_{[-\pi ,\pi ]^3}  \left( \frac{\partial f_N}{\partial t} -
  Q^R_B(f_N)\right) e^{-ik \cdot v} dv = 0.
\nonumber 
\end{gather}
After some computation we obtain
\begin{equation}
  {\hat Q}_k :=
  \int_{[-\pi,\pi]^3}
  \Q^R_B(f_N)
  e^{-i k \cdot v}\,dv
  = \sum_{\substack{l,m=-N/2\\l+m=k}}^{N/2} \f_l\,\f_m
  \hat{\beta}(l,m),\quad k=-N,\ldots,N,
  \label{eq:CF1}
\end{equation}
where $\hat{\beta }(l,m)=\hat{B}(l,m)-\hat{B}(m,m)$ are given by
\[
\hat{B}(l,m) = \int_{B_0(2\lambda\pi)}\int_{S^{2}} 
|q| \sigma(|q|, \cos\theta) e^{-i(l\cdot q^++m\cdot q^-)}\,d\omega\,dq.
\]
with 
\[
q^{+} = \frac12(q+\vert q\vert \omega), \quad
q^{-} = \frac12(q-\vert q\vert \omega).
\]
Finally, the set of ODEs is obtained:
\begin{gather}
  \frac{\partial \f_k}{\partial t}
  =
  \sum_{\substack{l,m=-N/2\nonumber\\l+m=k}}^{N/2} \hat{\beta }(l,m)\f_l\,\f_m
\nonumber 
\end{gather}
supplied with the initial condition
\begin{gather}
  \f_k(0) = {1\over (2\pi)^3 } \int_{[-\pi ,\pi ]^3} f_0(v)e^{-ik\cdot v} dv
\nonumber 
.
\end{gather}
Straightforward evaluation of (\ref{eq:CF1}) is expensive, especially
in three dimensions, the cost being of the order of $O(N^2)$. In order
to reduce the computational cost the so called Carleman representation
of (\ref{eq:boltzmann}) is used:
\[
  \Q_B (f,f)= \int_{\R^{3}} \int_{\R^{3}} {\tilde B}(x,y) 
  \delta(x \cdot y) 
  \left[ f(v + y) \, f(v+ x) - f(v+x+y) \, f(v) \right] \, dx \,
  dy,
\] 
with 
\[
  \tilde{B}(|x|,|y|) =
  2^{2} \, \sigma\left(\sqrt{|x|^2+|y|^2}, \frac{|x|}{\sqrt{|x|^2+|y|^2}} \right) \, (|x|^2+|y|^2)^{-\frac{1}{2}}.
\]
Under the $k$-th power inter particle force assumption (\ref{eq:J9hb})
this becomes
\begin{gather}
  \tilde{B}(|x|,|y|) = 4 C_\alpha (|x|^2 + |y|^2)^{-\frac{1-\alpha}{2} } 
  \label{eq:Btilde_pl}.
\end{gather}
The new quadrature formula is obtained:
\[
  \hat{Q}_k  =
  \sum_{\underset{l+m=k}{l,m=-N/2}}^{N/2} {\hat{\beta}}_F(l,m) \, \hat{f}_l \, \hat{f}_m, \ \ \
  k=-N,\ldots,N
\]
where ${\hat{\beta}}_F(l,m)=\hat{B}_F(l,m)-\hat{B}_F(m,m)$ are now 
given by
\[
  \hat{B}_F(l,m) = \int_{B_0(R)} \int_{B_0(R)}
  \tilde{B}(x,y) \, \delta(x \cdot y) \, 
  e^{i (l \cdot x+ m \cdot y)} \, dx \, dy.
\]
The next step is to identify a convolution structure. The goal is to
approximate ${\hat{\beta}}_F(l,m)$ by a sum
\begin{gather}
  {\hat{\beta}}_F(l,m) \simeq \sum_{p=1} ^{A} \alpha_p (l) \alpha' _p (m)
\nonumber 
,
\end{gather}
where $A$ represents a finite number of possible collision
directions. This is a discrete sum of $A$ convolutions and as a
consequence the computational cost of the algorithm is of the order of
$O(AN\log N)$.

This convolution can be obtained under assumption that
$\tilde{B}(|x|,|y|)$ is separable:
\begin{gather}
  \tilde{B}(|x|,|y|) = a(|x|) b(|y|)
\nonumber 
.
\end{gather}
This is the case is $\alpha $ is set to one in (\ref{eq:Btilde_pl}):
$\tilde{B}(|x|,|y|) = 4 C_\alpha $. In particular, when
$\tilde{B}(|x|,|y|) = 1$, this corresponds to the hard spheres
model. In this framework the following quadrature formula for
$\hat{B}_F(l,m)$ is obtained:
\begin{gather}
  \hat{B}_F(l,m) = {\pi^{2} \over A_1 A_2} \sum_{p,q=0}^{A_1,A_2}
  \alpha _{p,q}(l) \alpha _{p,q}^\prime (m)
\nonumber 
,
\end{gather}
where
\begin{gather}
  \alpha_{p,q} (l) = \phi_{R} ^3 \left(l \cdot e_{(\theta_p,\varphi_q)} \right),
  \hspace{0.8cm} \alpha^\prime _{p,q} (m) = \psi_{R} ^3
  \left(\Pi_{e_{(\theta_p,\varphi_q)} ^\bot} (m) \right),\nonumber\\
  \phi_{R} ^3 (s) = \int_{-R} ^R \rho \, e^{i \rho s} \, d\rho,
  \hspace{0.8cm} \psi_{R} ^3 (s) = \int_0 ^\pi \sin \theta \, \phi_{R} ^3 (s \cos \theta) \, d\theta, 
\nonumber 
\end{gather}
and the discrete angles $\theta _p$ and $\varphi _q$ are defined by
\begin{gather}
  (\theta_p,\varphi_q) = \Big(\frac{p \, \pi}{A_1}, \frac{q \, \pi}{A_2} \Big).
\nonumber 
\end{gather}

%
%

\bibliographystyle{abbrv}
\bibliography{biblio}

\begin{thebibliography}{10}

\bibitem{Aristov}
V.~V. Aristov and S.~A. Zabelok.
\newblock A deterministic method for solving the {B}oltzmann equation with
  parallel computations.
\newblock {\em Comput. Math. Math. Phys.}, 42(3):425--437, 2002.

\bibitem{Baranger}
C.~Baranger, J.~Claudel, N.~H{\'e}rouard, and L.~Mieussens.
\newblock Locally refined discrete velocity grids for stationary rarefied flow
  simulations.
\newblock {\em Journal of Computational Physics}, 257, Part A:572 -- 593, 2014.

\bibitem{BGK}
P.~L. Bhatnagar, E.~P. Gross, and M.~Krook.
\newblock {A model for collision processes in gases. I. Small amplitude
  processes in charged and neutral one-component systems}.
\newblock {\em Phys. Rev.}, 94(3):511--525, 1954.

\bibitem{bird}
G.~A. Bird.
\newblock {\em {Molecular gas dynamics and the direct simulation of gas
  flows}}.
\newblock Oxford University Press, 2nd edition, 1994.

\bibitem{bobylev}
A.~V. Bobylev, A.~Palczewski, and J.~Schneider.
\newblock On approximation of the {B}oltzmann equation by discrete velocity
  models.
\newblock {\em C. R. Acad. Sci. Paris Ser. I Math.}, 320(5):639--644, 1995.

\bibitem{bobylev_rjasanow97}
A.~V. Bobylev and S.~Rjasanow.
\newblock Difference scheme for the {B}oltzmann equation based on the fast
  {F}ourier transform.
\newblock {\em Eur. J. Mech. B Fluids}, 16(2):293--306, 1997.

\bibitem{Cf}
R.~E. Caflisch.
\newblock Monte carlo and quasi-{m}onte {c}arlo methods.
\newblock {\em Acta numerica}, 7:1--49, 1998.

\bibitem{CPima}
R.~E. Caflisch and L.~Pareschi.
\newblock {Towards a hybrid Monte Carlo method for rarefied gas dynamics}.
\newblock In {\em Transport in Transition Regimes}, pages 57--73. Springer,
  2004.

\bibitem{canuto:88}
C.~Canuto, M.~Hussaini, A.~Quarteroni, and T.~A. Zang.
\newblock {\em {Spectral methods in fluid dynamics}}.
\newblock Springer Series in Computational Physics. Springer-Verlag, New York,
  1988.

\bibitem{cercignani}
C.~Cercignani.
\newblock {\em The {B}oltzmann equation and its applications}, volume~67 of
  {\em Applied Mathematical Sciences}.
\newblock Springer-Verlag, New York, 1988.

\bibitem{CrSon1}
N.~Crouseilles, M.~Mehrenberger, and E.~Sonnendr{\"u}cker.
\newblock Conservative semi-{L}agrangian schemes for {V}lasov equations.
\newblock {\em Journal of Computational Physics}, 229(6):1927--1953, 2010.

\bibitem{CrSon}
N.~Crouseilles, T.~Respaud, and E.~Sonnendr{\"u}cker.
\newblock A forward semi-lagrangian method for the numerical solution of the
  vlasov equation.
\newblock {\em Computer Physics Communications}, 180(10):1730--1745, 2009.

\bibitem{FKS}
G.~Dimarco and R.~Loub{\`e}re.
\newblock {Towards an ultra efficient kinetic scheme. Part I: Basics on the
  {BGK} equation}.
\newblock {\em Journal of Computational Physics}, 255:680--698, 2013.

\bibitem{FKS_HO}
G.~Dimarco and R.~Loub{\`e}re.
\newblock {Towards an ultra efficient kinetic scheme. Part II: The high order
  case}.
\newblock {\em Journal of Computational Physics}, 255:699--719, 2013.

\bibitem{FKS_GPU}
G.~Dimarco, R.~Loub{\`e}re, and J.~Narski.
\newblock {Towards an ultra efficient kinetic scheme. Part III:
  High-performance-computing}.
\newblock {\em Journal of Computational Physics}, 284:22--39, 2015.

\bibitem{FKS_Boltz}
G.~Dimarco, R.~Loub{\`e}re, J.~Narski, and T.~Rey.
\newblock An efficient numerical method for solving the boltzmann equation in
  multidimensions.
\newblock {\em Submitted, arXiv:1608.08009}, 2016.

\bibitem{FKS_DD}
G.~Dimarco, R.~Loub{\`e}re, and V.~Rispoli.
\newblock {A multiscale fast semi-Lagrangian method for rarefied gas dynamics}.
\newblock {\em Journal of Computational Physics}, 291:99--119, 2015.

\bibitem{Dimarco_stiff1}
G.~Dimarco and L.~Pareschi.
\newblock High order asymptotic-preserving schemes for the {B}oltzmann
  equation.
\newblock {\em C. R. Math. Acad. Sci. Paris}, 350(9-10):481--486, 2012.

\bibitem{Dimarco_stiff2}
G.~Dimarco and L.~Pareschi.
\newblock Asymptotic preserving implicit-explicit {R}unge-{K}utta methods for
  nonlinear kinetic equations.
\newblock {\em SIAM J. Numer. Anal.}, 51(2):1064--1087, 2013.

\bibitem{DP-ACTA14}
G.~Dimarco and L.~Pareschi.
\newblock Numerical methods for kinetic equations.
\newblock {\em Acta Numer.}, 23:369--520, 2014.

\bibitem{filbet:2011Conv}
F.~Filbet and C.~Mouhot.
\newblock {{A}nalysis of spectral methods for the homogeneous {B}oltzmann
  Equation}.
\newblock {\em Trans. Amer. Math. Soc.}, 363:1947--1980, 2011.

\bibitem{FiMoPa:2006}
F.~Filbet, C.~Mouhot, and L.~Pareschi.
\newblock {Solving the Boltzmann equation in N log2 N}.
\newblock {\em SIAM J. Sci. Comput.}, 28(3):1029--1053, 2007.

\bibitem{Filbet2}
F.~Filbet and G.~Russo.
\newblock High order numerical methods for the space non-homogeneous
  {B}oltzmann equation.
\newblock {\em J. Comput. Phys.}, 186(2):457--480, Apr. 2003.

\bibitem{FilbetRusso}
F.~Filbet and G.~Russo.
\newblock Accurate numerical methods for the {B}oltzmann equation.
\newblock In {\em Modeling and computational methods for kinetic equations},
  pages 117--145. Springer, 2004.

\bibitem{Filbet}
F.~Filbet, E.~Sonnendr\"{u}cker, and P.~Bertrand.
\newblock Conservative numerical schemes for the {V}lasov equation.
\newblock {\em J. Comput. Phys.}, 172(1):166--187, Sept. 2001.

\bibitem{Frezzotti}
A.~Frezzotti, G.~P. Ghiroldi, and L.~Gibelli.
\newblock Direct solution of the boltzmann equation for a binary mixture on
  gpus.
\newblock {\em AIP Conference Proceedings}, 1333(1):884--889, 2011.

\bibitem{Frezzotti2}
A.~Frezzotti, G.~P. Ghiroldi, and L.~Gibelli.
\newblock Solving model kinetic equations on {GPU}s.
\newblock {\em Comput. \& Fluids}, 50:136--146, 2011.

\bibitem{Frezzotti3}
A.~Frezzotti, G.~P. Ghiroldi, and L.~Gibelli.
\newblock Solving the {B}oltzmann equation on {GPU}s.
\newblock {\em Comput. Phys. Comm.}, 182(12):2445--2453, 2011.

\bibitem{GambaL}
I.~M. Gamba and J.~R. Haack.
\newblock A conservative spectral method for the {B}oltzmann equation with
  anisotropic scattering and the grazing collisions limit.
\newblock {\em J. Comput. Phys.}, 270:40--57, 2014.

\bibitem{gamba}
I.~M. Gamba and S.~H. Tharkabhushanam.
\newblock {Spectral-Lagrangian methods for collisional models of
  non-equilibrium statistical states}.
\newblock {\em J. Comput. Phys.}, 228(6):2012--2036, Apr. 2009.

\bibitem{gamba:2010}
I.~M. Gamba and S.~H. Tharkabhushanam.
\newblock Shock and boundary structure formation by spectral-{L}agrangian
  methods for the inhomogeneous {B}oltzmann transport equation.
\newblock {\em J. Comput. Math.}, 28(4):430--460, 2010.

\bibitem{Gu}
Y.~G{\"u}{\c{c}}l{\"u} and W.~N.~G. Hitchon.
\newblock A high order cell-centered semi-{L}agrangian scheme for
  multi-dimensional kinetic simulations of neutral gas flows.
\newblock {\em J. Comput. Phys.}, 231(8):3289--3316, 2012.

\bibitem{Haack}
J.~Haack and I.~M. Gamba.
\newblock High performance computing with a conservative spectral boltzmann
  solver.
\newblock {\em 28th International Symposium on Rarefied Gas Dynamics 2012},
  1501:334--341, 2012.

\bibitem{Malkov}
E.~A. Malkov and M.~S. Ivanov.
\newblock Parallelization of algorithms for solving the boltzmann equation for
  gpu-based computations.
\newblock {\em AIP Conference Proceedings}, 1333(1):946--951, 2011.

\bibitem{Mieussens}
L.~Mieussens.
\newblock Discrete velocity model and implicit scheme for the {BGK} equation of
  rarefied gas dynamics.
\newblock {\em Mathematical Models and Methods in Applied Sciences},
  10(08):1121--1149, 2000.

\bibitem{MoPa:2006}
C.~Mouhot and L.~Pareschi.
\newblock Fast algorithms for computing the {B}oltzmann collision operator.
\newblock {\em Math. Comp.}, 75(256):1833--1852 (electronic), 2006.

\bibitem{Nanbu80}
K.~Nanbu.
\newblock Direct simulation scheme derived from the boltzmann equation. i.
  monocomponent gases.
\newblock {\em Journal of the Physical Society of Japan}, 49(5):2042--2049,
  1980.

\bibitem{Pal1}
A.~Palczewski and J.~Schneider.
\newblock Existence, stability, and convergence of solutions of discrete
  velocity models to the boltzmann equation.
\newblock {\em Journal of statistical physics}, 91(1-2):307--326, 1998.

\bibitem{pareschi:1996}
L.~Pareschi and B.~Perthame.
\newblock A fourier spectral method for homogeneous {B}oltzmann equations.
\newblock {\em Transport Theory Statist. Phys.}, 25(3):369--382, 1996.

\bibitem{PR2}
L.~Pareschi and G.~Russo.
\newblock Numerical solution of the {B}oltzmann equation {I: S}pectrally
  accurate approximation of the collision operator.
\newblock {\em SIAM J. Numer. Anal.}, 37(4):1217--1245, 2000.

\bibitem{Villani}
L.~Pareschi, G.~Toscani, and C.~Villani.
\newblock Spectral methods for the non cut-off {B}oltzmann equation and
  numerical grazing collision limit.
\newblock {\em Numer. Math.}, 93(3):527--548, 2003.

\bibitem{Shoucri}
M.~Shoucri and G.~Knorr.
\newblock Numerical integration of the {V}lasov equation.
\newblock {\em J. Computational Phys.}, 14(1):84--92, 1974.

\bibitem{titarev2012}
V.~Titarev.
\newblock Efficient deterministic modelling of three-dimensional rarefied gas
  flows.
\newblock {\em Communications in Computational Physics}, 12(01):162--192, 2012.

\bibitem{titarev2014construction}
V.~Titarev, M.~Dumbser, and S.~Utyuzhnikov.
\newblock Construction and comparison of parallel implicit kinetic solvers in
  three spatial dimensions.
\newblock {\em Journal of Computational Physics}, 256:17--33, 2014.

\bibitem{toro-book}
E.~F. Toro.
\newblock {\em Riemann solvers and numerical methods for fluid dynamics: a
  practical introduction}.
\newblock Springer Science \& Business Media, 2013.

\bibitem{wu2015influence}
L.~Wu, H.~Liu, Y.~Zhang, and J.~M. Reese.
\newblock Influence of intermolecular potentials on rarefied gas flows: Fast
  spectral solutions of the {B}oltzmann equation.
\newblock {\em Physics of Fluids (1994-present)}, 27(8):082002, 2015.

\bibitem{wu2013deterministic}
L.~Wu, C.~White, T.~J. Scanlon, J.~M. Reese, and Y.~Zhang.
\newblock {Deterministic numerical solutions of the Boltzmann equation using
  the fast spectral method}.
\newblock {\em Journal of Computational Physics}, 250:27--52, 2013.

\bibitem{wu2015kinetic}
L.~Wu, C.~White, T.~J. Scanlon, J.~M. Reese, and Y.~Zhang.
\newblock {A kinetic model of the Boltzmann equation for non-vibrating
  polyatomic gases}.
\newblock {\em Journal of Fluid Mechanics}, 763:24--50, 2015.

\bibitem{wu2015fast}
L.~Wu, J.~Zhang, J.~M. Reese, and Y.~Zhang.
\newblock {A fast spectral method for the Boltzmann equation for monatomic gas
  mixtures}.
\newblock {\em Journal of Computational Physics}, 298:602--621, 2015.

\end{thebibliography}

%
%

%

\end{document}